\documentclass[leqno, notitlepage]{amsart}

\usepackage{amsmath,amscd,amsfonts}
\usepackage[mathscr]{eucal}
\usepackage{epsfig}

\newcommand{\qrep}[1][d]{\ensuremath{\mathfrak{L}(#1)}}

\newcommand{\qrepc}[2]{\ensuremath{\mathfrak{L}(#1,#2)}}

\newcommand{\qrept}[1][d]{\ensuremath{\mathfrak{M}(#1)}}
\newcommand{\qreptc}[2]{\ensuremath{\mathfrak{M}(#1,#2)}}
\newcommand{\qreptr}[2][r]{\ensuremath{\mathfrak{M}^{#1}(#2)}}
\newcommand{\qreptcr}[3][r]{\ensuremath{\mathfrak{M}^{#1}(#2,#3)}}
\newcommand{\dqrept}[3]{\ensuremath{\mathfrak{M}(#1,#2,#3)}}
\newcommand{\dqreptr}[4][r]{\ensuremath{\mathfrak{M}^{#1}(#2,#3,#4)}}
\newcommand{\tv}[1][\mathbf{d}]{\ensuremath{\mathfrak{T}(#1)}}
\newcommand{\tvo}[1][\mathbf{d}]{\ensuremath{\mathfrak{T}_0(#1)}}
\newcommand{\B}{\ensuremath{\mathcal{B}}}
\newcommand{\id}{\ensuremath{\mathbf {1}}}
\newcommand{\twine}[2][\mathbf{d}_1,\ldots,\mathbf{d}_k]{\ensuremath{H_{#1}^{#2}}}
\newcommand{\cm}[2][\mathbf{d}_1,\ldots,\mathbf{d}_k]{\ensuremath{CM_{#1}^{#2}}}
\newcommand{\ocm}[2][\mathbf{d}_1,\ldots,\mathbf{d}_k]{\ensuremath{OCM_{#1}^{#2}}}
\newcommand{\lcm}[1][\mathbf{d}_1,\ldots,\mathbf{d}_k]{\ensuremath{LCM_{#1}}}
\newcommand{\olcm}[1][\mathbf{d}_1,\ldots,\mathbf{d}_k]{\ensuremath{OLCM_{#1}}}
\newcommand{\A}{\ensuremath{A_{\mathbf{w},\mathbf{r},\mathbf{n}}}}
\newcommand{\F}{\ensuremath{\id_{\A}}}
\newcommand{\f}[1][\mathbf{w}]{\ensuremath{f_{{#1},\mathbf{r},\mathbf{n}}}}
\newcommand{\K}[1][\mathbf{w}]{\ensuremath{k_{{#1},\mathbf{r},\mathbf{n}}}}
\newcommand{\gr}[2]{\ensuremath{Gr_{#1}^{#2}}}
\newcommand{\Ulong}{\ensuremath{U_q(\mathfrak{sl}_2)}}
\newcommand{\U}{\ensuremath{\mathbf{U}_q}}
\newcommand{\qchoose}[2]{\ensuremath{\begin{bmatrix} {#1} \\ {#2}
    \end{bmatrix}}}
\newcommand{\sv}[2][\mathbf{d}]{{\ensuremath{^{#1}v_{#2}}}}
\newcommand{\dv}[2][\mathbf{d}]{{\ensuremath{^{#1}v^{#2}}}}

\newcommand{\Perv}{\ensuremath{\mathcal{P}}}

\newcommand{\D}{\ensuremath{\mathcal{D}}}
\newcommand{\C}{\ensuremath{{\mathbb {C}}}} 
\newcommand{\Fq}{\ensuremath{{\mathbb {F}}_{q^2}}} 
\newcommand{\barFq}{\ensuremath{{\mathbb {{\bar F}}}_{q^2}}} 
\newcommand{\Z}{\mathbb {Z}} 
\newcommand{\N}{\mathbb {N}} 
\newcommand{\Proj}{\mathbb {P\,}} 
\newcommand{\s}[1]{|#1|}
\newcommand{\pair}[2]{\left<#1,#2\right>}
\newcommand{\Span}{\operatorname{Span}}
\newcommand{\End}{\operatorname{End}}
\newcommand{\Hom}{\operatorname{Hom}}

\newtheorem{lemma}{Lemma}[subsection]
\newtheorem{prop}{Proposition}[subsection]
\newtheorem{theorem}{Theorem}[subsection]
\newtheorem{conjecture}{Conjecture}[subsection]

\DeclareMathOperator{\rank}{rank}
\DeclareMathOperator{\im}{im} 

\allowdisplaybreaks
\begin{document}
\title{The Tensor Product of Representations of \Ulong\ Via Quivers}
\author{Alistair Savage}
\address{ Department of Mathematics, Yale University,
P.O. Box 208283, New Haven, CT 06520-8283, USA }
\email{alistair.savage@yale.edu}
\date{\today}

\begin{abstract}

Using the tensor product variety introduced in \cite{malkin} and
\cite{nakajima:01}, the complete structure of the tensor product of a
finite number of integrable highest weight modules of
\Ulong\ is recovered.  In particular, the
elementary basis, Lusztig's canonical basis, and the basis
adapted to the decomposition of the tensor product
into simple modules are all exhibited as distinguished elements of
certain spaces of
invariant functions on the tensor product variety.  For the latter
two bases, these distinguished elements are closely related to the
irreducible components of the tensor product variety.  The space of
intertwiners is also interpreted geometrically.

\end{abstract}

\maketitle
\tableofcontents

\section*{Introduction}
The purpose of this paper is to obtain a geometric description of the
tensor product of a
finite number of integrable highest weight representations
of \Ulong\ using quiver varieties.  The definition of a
\emph{tensor product variety}
corresponding to the tensor product of a finite number of
integrable highest weight representations of a Lie algebra
$\mathfrak{g}$ of ADE type was
introduced in \cite{malkin} and
\cite{nakajima:01} (see also \cite{varagnolo-vasserot:01} for a
geometric description of the tensor product).  There it is
demonstrated that the set of irreducible components of the tensor
product variety can be
equipped with the structure of a $\mathfrak{g}$-crystal isomorphic to
the crystal of the canonical basis in the tensor product representation.

In this paper, we consider the specific case $\mathfrak{g} =
\mathfrak{sl}_2$ and recover the entire structure (as opposed to the
crystal structure alone)
of \Ulong\ via the tensor product variety.  Our definition of the
tensor product variety differs
slightly from that of \cite{malkin} and \cite{nakajima:01} in that we
consider our varieties over the finite field \Fq\ with $q^2$ elements
(or its algebraic
closure \barFq) rather than over \C.  The reader who is only interested
in representations of $sl_2$, rather than its associated quantum
group, may replace \Fq\ by \C\ and set $q=1$ everywhere.  With a few
obvious modifications, the arguments of the paper still hold.
Let $\mathbf{d} \in (\Z_{\ge
  0})^k$.  We find three distinct spaces, $\mathcal{T}_0(\mathbf{d})$,
$\mathcal{T}_c(\mathbf{d})$, and
$\mathcal{T}_s(\mathbf{d})$, of invariant (with respect to a
natural group action)
functions on the tensor
product variety
\tv, each
isomorphic to $V_{\mathbf{d}_1} \otimes \dots \otimes
  V_{\mathbf{d}_k}$.  In each space we define a natural basis.  These
three bases, $\B_e$, $\B_c$, and $\B_s$, correspond respectively to
the elementary basis, Lusztig's canonical basis \cite{lusztig:92}, and a basis
compatible with the
decomposition of $V_{\mathbf{d}_1} \otimes \dots
\otimes V_{\mathbf{d}_k}$ into a direct sum of irreducible modules.
The two bases $\B_c$ and $\B_s$ are characterized by their relation to the
irreducible components of \tv.  We define the irreducible
components of \tv\ (defined over \Fq) to be the \Fq\ points of the irreducible components of $\tv'$ (the corresponding variety defined over \barFq).
We then define the \emph{dense
  points} of an irreducible component of \tv\ to be the \Fq\ points of a
certain dense subset of the corresponding irreducible component of
$\tv'$.  Distinct
elements of the basis $\B_c$ and $\B_s$ are supported on distinct irreducible
components of \tv\ and equal to a non-zero constant on the set of
dense points of that
irreducible component (see
Theorems \ref{thm:B_c} and \ref{thm:char-B_s}).  However, the supports
of the elements of $\B_s$ are disjoint whereas the supports of the
elements of $\B_c$ are not.
We also find a geometric
description of the space of intertwiners $\twine{\mu} =
\Hom_{\Ulong} (V_{\mathbf{d}_1} \otimes \dots \otimes
V_{\mathbf{d}_k}, V_\mu)$.  A natural basis $\B_I$ of this space is
again characterized by its relation
to the irreducible components of
\tv.

An important tool used in the development and proof of the
results of this paper is the graphical
calculus of intertwiners of \Ulong\ introduced by Penrose,
Kauffman and
others.  This graphical calculus is expanded
in \cite{frenkel-khovanov:97} and used to prove various results
concerning Lusztig's canonical basis.
The present paper can be considered a ``geometrization'' of these results.

In Section~\ref{sec:perv-conj} we conjecture a characterization of the
basis $\B_c$ as the image of certain intersection cohomology
sheaves of \tv\ under a particular functor from the space of
constructible
 semisimple perverse sheaves on \tv\ to the space of invariant
functions on \tv.  Since the definition of
$\mathcal{T}_c(\mathbf{d})$ relies on the graphical calculus of
intertwiners of \Ulong\ (and no such graphical calculus exists for
more general Lie algebras), this conjecture should play a key role in
the possible extension of the results of this paper to a more general set of
Lie algebras (for instance, those of type ADE).

This paper is organized as follows.  Section~\ref{sec:U_reps_review}
contains a review of
\Ulong\ and its representations, Nakajima's quiver varieties, and the
graphical calculus of intertwiners of \Ulong.  The tensor product variety is
defined in Section~\ref{sec:geom_tensor} where the spaces
$\mathcal{T}_0(\mathbf{d})$
and $\mathcal{T}_c(\mathbf{d})$ are introduced, an isomorphism between
the two is given, and various results
concerning these spaces and their distinguished bases $\B_e$ and
$\B_c$ are proved.
Section~\ref{sec:intertwiners} is concerned with a geometric
realization of the space of
intertwiners and the decomposition of the tensor product
representation into a direct sum of irreducible modules (via the space
$\mathcal{T}_s(\mathbf{d})$ and the distinguished basis $\B_s$).  It
is concluded with the discussion of an isomorphism between the spaces
$\mathcal{T}_c(\mathbf{d})$ and $\mathcal{T}_s(\mathbf{d})$.

The notation used in the description of quiver varities is not
standardized.  Lusztig denotes the fixed vector space by $D$ and the
subpace by $V$ while Nakajima denotes these objects by $W$ and $V$
respectively.  Since we wish to use the notation $V_n$ for certain
\Ulong\ modules (to agree with the notation of
\cite{frenkel-khovanov:97}), we denote the fixed vector space by $D$
and the subspace by $W$.  We hope that this will not cause confusion among
those readers familiar with the work of Lusztig and Nakajima.

Throughout this paper the topology is the
Zariski topology and the ground field is $\barFq$ unless otherwise
specified.  However, we will usually deal with varieties defined over \Fq\
and consider the corresponding set of \Fq -rational points.
Thus, for instance, $\Proj^n = \Proj^n \Fq$ and a vector
space is an \Fq\ vector space.  A
function on an algebraic variety is a function into $\C(q)$, the field
of rational functions in an indeterminate $q$.  The span of a set of
such functions is their $\C(q)$-span. The support of a
function $f$ is defined to be
the set $\{x \, |\, f(x) \ne 0\}$ and \emph{not} the closure of this set.

I would like to thank I. Frenkel for suggesting the topic of this
paper and for his help during its development.  I am also grateful to
A. Malkin, O. Schiffmann, M. Khovanov and H. Nakajima for very helpful
discussions.
This work was
supported in part by the Natural Science and Engineering Research Council of
Canada.


\section{The Quantum Group \Ulong\ and its Representations}
\label{sec:U_reps_review}
\subsection{The Hopf Algebra Structure of \Ulong}
Let $\C(q)$ be the field of rational functions in an indeterminate $q$
and define $\bar{\ } : \C(q) \rightarrow \C(q)$ to be the $\C$-algebra
involution such that $\overline{q^n} = q^{-n}$ for all $n$.  The
quantum group \Ulong\ (which we will denote by \U) is the associative
algebra over $\C(q)$ with
generators $E,F,K,K^{-1}$ and relations
\begin{align*}
KK^{-1} &= K^{-1}K \\
KE &= q^2EK \\
KF &= q^{-2}FK \\
EF - FE &= \frac{K - K^{-1}}{q-q^{-1}}.
\end{align*}
The comultiplication and counit of the Hopf algebra structure of \U\
are given by
\begin{align*}
&\Delta K^{\pm 1} = K^{\pm 1} \otimes K^{\pm 1} \\
&\Delta E = E \otimes 1 + K \otimes E \\
&\Delta F = F \otimes K^{-1} + 1 \otimes F \\
\end{align*}
and
\begin{align*}
&\eta(K^{\pm 1}) = 1 \\
&\eta(E) = \eta(F) = 0
\end{align*}
respectively.  Although an explicit expression for the antipode
exists, we will not need it in this paper.

Let us introduce two involutions of \U.  The first one is the \emph{Cartan
involution}, denoted by $\omega$, which acts as follows:
\begin{align*}
&\omega(E) = F,\quad \omega(F) = E, \quad \omega(K^{\pm 1}) = K^{\pm
  1},\quad \omega(q^{\pm 1}) = q^{\pm 1} \\
&\omega(xy) = \omega(y)\omega(x),\quad x,y \in \U.
\end{align*}
The second, denoted by $\sigma$, is called the \emph{``bar''
  involution} and is defined by
\begin{align*}
&\sigma(E) = E,\quad \sigma(F) = F,\quad \sigma(K^{\pm 1}) = K^{\mp
  1},\quad \sigma(q^{\pm 1}) = q^{\mp 1} \\
&\sigma(xy) = \sigma(x) \sigma (y),\quad x,y \in \U.
\end{align*}
Using $\sigma$ we can define a second comultiplication $\bar \Delta$
by
\begin{equation*}
\bar\Delta (x) = (\sigma \otimes \sigma)\Delta(\sigma(x)),\quad x\in\U
\end{equation*}
which implies
\begin{align*}
&\bar\Delta K^{\pm 1} = K^{\pm 1} \otimes K^{\pm 1} \\
&\bar\Delta E = E \otimes 1 + K^{-1} \otimes E \\
&\bar\Delta F = F \otimes K + 1 \otimes F.
\end{align*}


\subsection{Irreducible Representations of \Ulong}
\label{subsec:reps}
Any finite dimensional irreducible \U-module $V$ is generated by a highest
weight vector, $v$, of
weight $\varepsilon q^d$ where $\varepsilon = \pm 1$ and $d = \dim (V)
-1$ \cite{kassel}.  In this paper we consider those representations with
$\varepsilon = +1$.  Let $v_{d-2k} = F^kv/[k]!$ where
\begin{align*}
&[k] = (q^k -
q^{-k})/(q-q^{-1}) = q^{-k+1} + q^{-k+3} + \dots + q^{k-1}, \\
&[k]! = [1][2]\cdots[k].
\end{align*}
Then $v_{d-2k} = 0$ for $k > d$ and $\{v=v_d,
v_{d-2}, \dots, v_{-d}\}$ is a basis of $V$.  We denote this
representation by $V_d$.  The action of \U\ on $V_d$ is given by
\begin{equation}
\label{U_action_V_d}
\begin{split}
& K^{\pm 1} v_m = q^{\pm m} v_m \\
& Ev_m = \left[ \frac{d+m}{2} + 1 \right] v_{m+2} \\
& Fv_m = \left[ \frac{d-m}{2} + 1 \right] v_{m-2}.
\end{split}
\end{equation}

Define a bilinear symmetric pairing on $V_d$ by requiring
\begin{equation*}
\pair{xu}{v} = \pair{u}{\omega(x)v},\quad \pair{v_d}{v_d} = 1,\quad
u,v \in V_d \text{ and } x \in \U.
\end{equation*}
It follows that
\begin{equation*}
\pair{v_{d-2k}}{v_{d-2l}} = \delta_{k,l} \qchoose{d}{k}
\end{equation*}
where
\begin{equation*}
\qchoose{d}{k} = \frac{[d]!}{[k]![d-k]!}.
\end{equation*}

Let $\{v^{d-2k}\}_{k=0}^d$ be the basis dual to $\{v_{d-2k}\}_{k=0}^d$ with
respect to the form $\pair{\ }{\ }$.  Then
\begin{equation*}
v^{d-2k} = \qchoose{d}{k}^{-1} v_{d-2k}
\end{equation*}
and the action of \U\ in the dual basis is
\begin{align*}
&K^{\pm 1}v_m = q^{\pm m} v_m \\
&E v^m = \left[ \frac{d-m}{2} \right] v^{m+2} \\
&F v^m = \left[ \frac{d+m}{2} \right] v^{m-2}.
\end{align*}


\subsection{Geometric Realization of Irreducible Representations of \Ulong}
\label{sec:geom_reps}
We recall here Nakajima's quiver variety construction of finite
dimensional irreducible
representations of Kac-Moody algebras associated to symmetric Cartan matrices
\cite{nakajima:94,nakajima:98} in the specific
case of \Ulong.  In order to introduce the quantum parameter $q$, some
of our definitions differ slightly from those in \cite{nakajima:94,
  nakajima:98}.  Since the Dynkin diagram of $\mathfrak{sl}_2$
consists of a single vertex and no edges, the definition of the quiver variety
simplifies considerably.
Fix vector spaces $W$ and $D$ of dimensions $w$ and $d$ respectively
and consider the
variety
\[
\mathbf{M}(w,d) = \Hom (D,W) \oplus \Hom (W,D).
\]
The two components of an element of $\mathbf{M}(w,d)$ will be denoted by $f_1$
and $f_2$ respectively.  $GL(W)$ acts on $\mathbf{M}(w,d)$ by
\[
(f_1,f_2) \mapsto g(f_1,f_2) \stackrel{\text{def}}{=} (gf_1,
f_2g^{-1}),\, g \in GL(W).
\]
Define the map $\mu: \mathbf{M}(w,d) \to \End W$ by
\[
\mu(f_1,f_2) = f_1 f_2.
\]
Let $\mu^{-1}(0)$ be the algebraic variety defined as the zero set of
$\mu$.  We say a point $(f_1, f_2)$ of $\mu^{-1}(0)$ is \emph{stable}
if $f_2$ is injective.  The \emph{quiver variety} is then given by
\[
\{(f_1,f_2) \in \mu^{-1}(0)\, |\, (f_1,f_2) \text{ is stable}\}/GL(W).
\]
Via the map $(f_1,f_2) \mapsto (\im f_2, f_2 f_1)$, this variety is
seen to be isomorphic to the variety
\[
\qreptc{w}{d}=\{(W,t)\, |\, W\subset D,\, \dim W = w,\, t\in \End D,\, \im t
\subset W  \subset \ker t\}.
\]
Note that the condition $\im t
\subset W \subset \ker t$ implies $t^2 = 0$.
Let
\[
\qrept= \bigcup_{w} \qreptc{w}{d} = \{(W,t)\, |\, W\subset D,\, t\in
\End D,\, \im t \subset W  \subset \ker t\}.
\]
and
\begin{multline*}
\dqrept{w}{w+1}{d} = \{(U,W,t)\, |\, t\in \End D,\, \im t \subset U
\subset W \subset \ker t,\, \\
\dim U=w,\, \dim W=w+1\}.
\end{multline*}
We then have the projections
\[
\qrept \stackrel{\pi_1}{\longleftarrow} \bigcup_w \dqrept{w}{w+1}{d}
\stackrel{\pi_2}{\longrightarrow} \qrept
\]
given by $\pi_1(U,W,t)=(U,t)$ and $\pi_2(U,W,t)=(W,t)$.

For a subset $Y$ of a variety $A$, let $\id_Y$ denote the
function on $A$ which takes the value 1 on $Y$ and 0 elsewhere.
Note that since our varieties are defined over \Fq, they consist of a
finite number of (\Fq - rational) points.  Let $\chi_q(Y)$ denote the Euler
characteristic of the algebraic variety Y, which is merely the number of points
in $Y$.
For a map $\pi$ between algebraic varieties $A$ and $B$, let $\pi_!$
\cite{macpherson:74} denote the
map between the abelian groups of functions on $A$ and
$B$ given by
\begin{align*}
\pi_!(f)(x) &= \sum_{y \in \pi^{-1}(x)} f(y) \\
\Rightarrow \pi_!(\id_Y)(x) &= \chi_q(\pi^{-1}(x)\cap Y),\; Y \subset A
\end{align*}
and let $\pi^*$ be the pullback map from
functions on $B$ to functions on $A$ acting as $\pi^*f(x)=f(\pi(x))$.

We then define the action of $E$, $F$ and $K^{\pm 1}$ on the set of
functions on \qrept\ by
\begin{equation}
\label{EF_action}
\begin{split}
&Ef = q^{-\dim(\pi_1^{-1}(\cdot))} (\pi_1)_!\pi_2^*f \\
&Ff = q^{-\dim(\pi_2^{-1}(\cdot))} (\pi_2)_!\pi_1^*f \\
&K^{\pm 1}f = q^{\pm (d-2\dim (\cdot))} f
\end{split}
\end{equation}
where the notation means that for a function $f$ on
\qrept\ and $(W,t) \in \qrept$,
\begin{equation}
\begin{split}
&Ef(W,t) = q^{-\dim(\pi_1^{-1}(W,t))} (\pi_1)_!\pi_2^*f(W,t) \\
&Ff(W,t) = q^{-\dim(\pi_2^{-1}(W,t))} (\pi_2)_!\pi_1^*f(W,t) \\
&K^{\pm 1}f(W,t) = q^{\pm (d-2\dim W)} f(W,t).
\end{split}
\end{equation}

Let
\begin{align*}
&\qreptr{d} = \{(W,t) \in \qrept\, |\, \rank t=r\} \\
&\qreptcr{w}{d} = \{(W,t) \in \qreptc{w}{d}\, |\, \rank t = r\} \\
&\mathcal{M}^r(w,d) = \C(q) \id_{\qreptcr{w}{d}} \\
&\mathcal{M}^r(d) = \bigoplus_{w} \mathcal{M}^r(w,d) \\
&\mathcal{M}(w,d) = \bigoplus_r \mathcal{M}^r(w,d) \\
&\mathcal{M}(d) = \bigoplus_w \mathcal{M}(w,d).
\end{align*}
Also, let us introduce the following notation for Grassmanians:
\[
\gr{w}{d} = \{W \subset (\Fq)^d\, |\, \dim W = w\} .
\]

\begin{prop}
\label{prop:qrep_isom}
The action of \U\ defined by \eqref{EF_action} endows
$\mathcal{M}^r(d)$ (and
hence $\mathcal{M}(d)$) with the structure of a \U-module and
the map $\id_{\qreptcr{w}{d}} \mapsto v_{d-2w}$ (extended by
linearity) is an isomorphism
$\mathcal{M}^r(d) \cong V_{d-2r}$ of \U -modules.
\end{prop}

To prove this proposition, we will need the following lemmas.

\begin{lemma}
\label{lem:factor_isom}
For vector spaces $W \subset D$, $\{U\, |\, W \subset U
\subset D,\, \dim U = u\} \cong \gr{u - \dim W}{\dim D - \dim W}$.
\end{lemma}
\begin{proof}
This follows immediately from the fact that
\[
\{U\, |\, W \subset U \subset D,\, \dim U = u\} \cong \{U'\, |\, U' \subset
  D/W,\, \dim U' = u - \dim W\}
\]
via the map $U \mapsto U' = U/W$.
\end{proof}

\begin{lemma}
\label{lem:chi}
$\chi_q(\Proj^n) = \sum_{i=0}^n q^{2i}$
\end{lemma}
\begin{proof}
This follows from simply counting the number of possible one
dimensional subspaces of $\Proj^n$.
\end{proof}

\begin{proof}[Proof of Proposition~\ref{prop:qrep_isom}]
If $(W,t) \in \qreptcr{w}{d}$ then
\begin{align*}
E\id_{\qreptcr{w+1}{d}}(W,t) &= q^{-\dim(\pi_1^{-1}(W,t))}
(\pi_1)_! \pi_2^* \id_{\qreptcr{w+1}{d}} (W,t) \\
&= q^{-\dim(\{U\, |\, W \subset U \subset \ker t,\, \dim U=w+1\})}
(\pi_1)_! \id_{\dqreptr{w}{w+1}{d}} (W,t) \\
&= q^{-\dim(\gr{1}{d-w-r})} \chi_q (\pi_1^{-1}(W,t) \cap
\dqreptr{w}{w+1}{d}) \\
&= q^{-\dim(\Proj^{d-w-r-1})} \chi_q (\{U\, |\, W \subset U \subset
\ker t,\, \dim U=w+1\}) \\
&= q^{-(d-w-r-1)} \chi_q (\gr{1}{d-w-r}) \\
&= q^{-(d-w-r-1)} \chi_q (\Proj^{d-w-r-1}) \\
&= q^{-(d-w-r-1)} \sum_{i=0}^{d-w-r-1} q^{2i} \\
&= q^{-(d-w-r-1)} + q^{-(d-w-r-1)+2} + \dots + q^{d-w-r-1} \\
&= [d-w-r]
\end{align*}
and $E\id_{\qreptcr{w+1}{d}}(W,t)=0$ otherwise.  So
$E\id_{\qreptcr{w+1}{d}} = [d-w-r]\id_{\qreptcr{w}{d}}$.  Similarly,
if $(W,t) \in \qreptcr{w+1}{d}$,
\begin{align*}
F\id_{\qreptcr{w}{d}}(W,t) &=
q^{-\dim(\pi_2^{-1}(W,t))} (\pi_2)_! \pi_1^* \id_{\qreptcr{w}{d}} (W,t) \\
&= q^{-\dim(\{U\, |\, \im t \subset U \subset W,\, \dim U = w\})}
(\pi_2)_! \id_{\dqreptr{w}{w+1}{d}} (W,t) \\
&= q^{-\dim(\gr{w-r}{w+1-r})} \chi_q (\pi_2^{-1}(W,t) \cap
\dqreptr{w}{w+1}{d}) \\
&= q^{-\dim(\Proj^{w-r})} \chi_q (\{U\, |\, \im t \subset U \subset
W,\, \dim U = w\}) \\
&= q^{-(w-r)} \chi_q (\Proj^{w-r}) \\
&= q^{-(w-r)} \sum_{i=0}^{w-r} q^{2i}\\
&= q^{-(w-r)} + q^{-(w-r)+2} + \dots + q^{w-r} \\
&= [w+1-r]
\end{align*}
and $F\id_{\qreptcr{w}{d}}(W,t)=0$ otherwise.  So
$F\id_{\qreptcr{w}{d}} = [w+1-r]\id_{\qreptcr{w+1}{d}}$.
It is obvious that 
\begin{equation}
K^{\pm 1}\id_{\qreptcr{w}{d}} = q^{\pm(d-2w)} \id_{\qreptcr{w}{d}}.
\end{equation}

Now, $\qreptcr{w}{d}=\emptyset$ unless $r \le w \le d-r$ due to the
requirement $\im t \subset W \subset \ker t$ in the definition of
\qreptcr{w}{d}. Thus $\mathcal{M}^r(d)=\bigoplus_{w=r}^{w=d-r}
\mathcal{M}^r(w,d)$.

Comparing the above calculations to \eqref{U_action_V_d}, the result follows.
\end{proof}

So $\mathcal{M}(d)$ is isomorphic to the direct sum of the irreducible
representations of highest weight $d-2r$ where $0 \le r \le d/2$ since these
are the possible ranks of $t$ (recall that $t^2=0$).

Let $\qrep = \qreptr[0]{d}$.  Then \qrep\ is isomorphic to the algebraic
variety of all subspaces $W\subset D$, which is a union of
Grassmanians.  Let
\[
\qrepc{w}{d} = \qreptcr[0]{w}{d} = \{W \subset D\, |\, \dim W = w\}
\cong \gr{w}{d}.
\]
and
\[
\mathcal{L}(w,d) = \mathcal{M}^0(w,d)= \C(q) \id_{\qrepc{w}{d}},\quad
\mathcal{L}(d)= \mathcal{M}^0(d) = \bigoplus_{w=1}^d \mathcal{L}(w,d).
\]
We see from Proposition~\ref{prop:qrep_isom} that the action of \U\
defined by \eqref{EF_action} endows $\mathcal{L}(d)$ with the
structure of the
irreducible module $V_d$ via the isomorphism $\id_{\qrepc{w}{d}} \mapsto
v_{d-2w}$ (extended by linearity).  Note that for $(W,t) \in \qrept$,
we can think of $t$ as belonging to $\Hom(D/W,W)$ and thus
\qrept\ is the cotangent bundle of \qrep.


\subsection{Tensor Products and the Graphical Calculus of Intertwiners}
\label{sec:crossmatch}

We define the bilinear pairing of $V_{\mathbf{d}_1} \otimes \dots
\otimes V_{\mathbf{d}_k}$ with $V_{\mathbf{d}_k} \otimes \dots
\otimes V_{\mathbf{d}_1}$ by
\begin{equation*}
\pair{v_{i_1} \otimes \dots \otimes v_{i_k}}{v^{l_k} \otimes \dots
  \otimes v^{l_1}} = \delta_{i_1}^{l_1} \dots \delta_{i_k}^{l_k}.
\end{equation*}
Then
\begin{multline*}
\pair{\Delta^{n-1}(x) v_{i_1} \otimes \dots \otimes v_{i_k}}{v^{l_k}
  \otimes \dots \otimes v^{l_1}} \\
= \pair{v_{i_1}
  \otimes \dots \otimes v_{i_k}}{\bar\Delta^{n-1}(\omega(x)) v^{l_k}
  \otimes \dots \otimes v^{l_1}}.
\end{multline*}

Lusztig's canonical basis of the tensor product is described in
\cite{lusztig:92}.  We refer
the reader to this article or the overview in
\cite{frenkel-khovanov:97}, Section 1.5, for the definition of this basis.
As in \cite{frenkel-khovanov:97} and \cite{lusztig:92},\
we denote the elements of Lusztig's
canonical basis by $v_{i_1} \diamondsuit \cdots \diamondsuit v_{i_k}$ and their
dual by $v_{i_1} \heartsuit \cdots \heartsuit v_{i_k}$.  The dual is
defined with respect to the form $\pair{\ }{\ }$:
\begin{equation*}
\pair{v_{i_1} \diamondsuit \cdots \diamondsuit v_{i_k}}{v^{l_k} \heartsuit \cdots
  \heartsuit v^{l_1}} = \delta_{i_1}^{l_1} \dots \delta_{i_k}^{l_k}.
\end{equation*}
When we wish to make explicit to which representation a vector
belongs, we use the notation $\sv[d]{k},\, \dv[d]{k} \in V_d$

To simplify notation, we make the following definitions
\begin{align*}
\otimes \sv{\mathbf{w}} &= \sv[\mathbf{d}_1]{\mathbf{d}_1 -
  2\mathbf{w}_1} \otimes \dots \otimes \sv[\mathbf{d}_k]{\mathbf{d}_k -
  2\mathbf{w}_k} \\
\diamondsuit \sv{\mathbf{w}} &= \sv[\mathbf{d}_1]{\mathbf{d}_1 -
  2\mathbf{w}_1} \diamondsuit \dots \diamondsuit \sv[\mathbf{d}_k]{\mathbf{d}_k -
  2\mathbf{w}_k} \\
\otimes \dv{\mathbf{w}} &= \dv[\mathbf{d}_1]{\mathbf{d}_1 -
  2\mathbf{w}_1} \otimes \dots \otimes \dv[\mathbf{d}_k]{\mathbf{d}_k -
  2\mathbf{w}_k} \\
\heartsuit \dv{\mathbf{w}} &= \dv[\mathbf{d}_1]{\mathbf{d}_1 -
  2\mathbf{w}_1} \heartsuit \dots \heartsuit \dv[\mathbf{d}_k]{\mathbf{d}_k -
  2\mathbf{w}_k}
\end{align*}
where $\mathbf{d}, \mathbf{w} \in (\Z_{\ge 0})^k$.

We can extend the bar involution $\sigma$ to tensor products of
irreducible representations as
follows.  Define
\[
\sigma \left(f(q) \left(\otimes
    \sv{\mathbf{w}}\right)\right) = f(q^{-1}) \left(\otimes
    \sv{\mathbf{w}}\right)
\]
and extend by $\C$-linearity.  Then $\sigma$ is an isomorphism from
$V_{\mathbf{d}_1} \otimes \dots \otimes V_{\mathbf{d}_k}$ to itself and
\begin{equation}
\label{eq:sigma}
\sigma (\Delta^{(k-1)}(x) (v)) = ((\sigma \otimes \dots \otimes
\sigma)(\Delta^{(k-1)}x))(\sigma v)
\end{equation}
for $x \in \U$ and $v \in V_{\mathbf{d}_1} \otimes \dots \otimes
V_{\mathbf{d}_k}$.

We now recall some results on the graphical calculus of tensor
products and intertwiners.  For a more complete treatment, see
\cite{frenkel-khovanov:97}.
In the graphical calculus, $V_d$ is depicted by a box
marked $d$ with $d$ vertices.  To depict $\cm{\mathbf{a}_1, \dots,
  \mathbf{a}_l}$, we place the
boxes representing the
$V_{\mathbf{d}_i}$ on a horizontal line and the boxes representing the
$V_{\mathbf{a}_i}$ on another horizontal line lying above the first one.
$\cm{\mathbf{a}_1, \dots, \mathbf{a}_l}$ is then the set of non-intersecting curves
(up to isotopy) connecting the vertices of the boxes such that the following
conditions are satisfied:

\begin{enumerate}
\item Each curve connects exactly two vertices.
\item Each vertex is the endpoint of exactly one curve.
\item No curve joins a box to itself.
\item The curves lie inside the box bounded by the two horizontal
  lines and the vertical lines through the extreme right and left points.
\end{enumerate}

\begin{figure}
\begin{center}
\epsfig{file=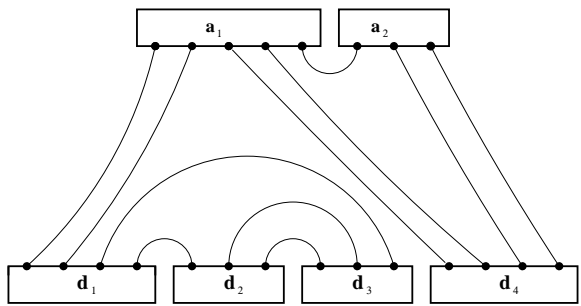,width=3in}
\caption{A crossingless match \label{fig:cm}}
\end{center}
\end{figure}

\begin{figure}
\begin{center}
\epsfig{file=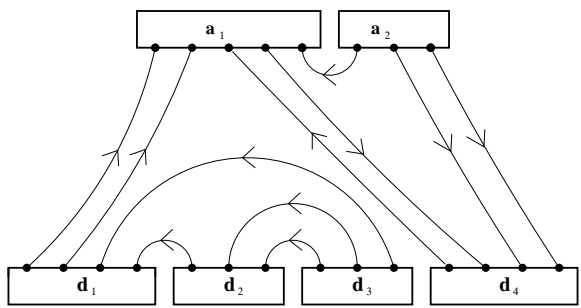,width=3in}
\caption{An oriented crossingless match \label{fig:ocm}}
\end{center}
\end{figure}

An example is given in Figure~\ref{fig:cm}.
We call the curves joining two lower boxes \emph{lower curves}, those
joining two upper boxes \emph{upper curves} and those joining a lower
and an upper box \emph{middle curves}.  We define the set of oriented
crossingless matches \ocm{\mathbf{a}_1, \dots, \mathbf{a}_l}\ to be the set of
elements of \cm{\mathbf{a}_1, \dots, \mathbf{a}_l}\ along with an orienation of the
curves such that all upper and lower curves are oriented to the left
and all middle curves are oriented
so that those oriented down are to the right of those oriented up.
See Figure~\ref{fig:ocm}.

As shown in \cite{frenkel-khovanov:97}, the set of crossingless
matches \cm{\mathbf{a}_1, \dots, \mathbf{a}_l}\ is in
one to one correspondence with a basis of the set of intertwiners
\[\twine{\mathbf{a}_1,
  \dots, \mathbf{a}_l} = \Hom_{\U} \left( V_{\mathbf{d}_1} \otimes \dots
\otimes V_{\mathbf{d}_k}, V_{\mathbf{a}_1} \otimes \dots \otimes V_{\mathbf{a}_l}
\right).
\]
The matrix coefficients of the intertwiner associated to a particular
crossingless match are given by Theorem~2.1 of \cite{frenkel-khovanov:97}.
Note that these are intertwiners in the dual basis and
thus commute with the action of \U\ on the tensor product given by ${\bar
  \Delta}^{(k-1)}$.  Let $\tilde \gamma$ be such an intertwiner
and define $\gamma = \sigma \tilde \gamma \sigma$.  Then for $x \in
\U$ and $v \in V_{\mathbf{d}_1} \otimes \dots \otimes V_{\mathbf{d}_k}$,
\begin{align*}
\gamma \Delta^{(k-1)} (x) (v) &= \sigma \tilde \gamma \sigma
\Delta^{(k-1)}(x) (v) \\
&= \sigma \tilde \gamma ((\sigma \otimes \dots \otimes \sigma) \Delta^{(k-1)}
(x)) (\sigma v) \\
&= \sigma \tilde \gamma {\bar \Delta^{(k-1)}} (\sigma x) (\sigma v) \\
&= \sigma {\bar \Delta^{(k-1)}} (\sigma x) \tilde \gamma (\sigma v) \\
&= \sigma ((\sigma \otimes \dots \otimes \sigma) \Delta^{(k-1)} (x))
\sigma \gamma (v) \\
&= \Delta^{(k-1)} (x) \gamma (v).
\end{align*}
Thus $\gamma$ is an intertwiner in the usual basis commuting with the
action of \U\ given by $\Delta^{(k-1)}$.

\begin{figure}
\begin{center}
\epsfig{file=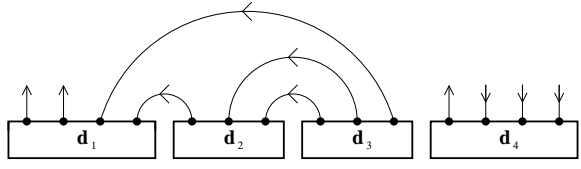,width=3in}
\caption{An oriented lower crossingless match \label{fig:olcm}}
\end{center}
\end{figure}

We will also need to define the set of \emph{lower crossingless matches}
\lcm\ and \emph{oriented lower crossingless matches} \olcm.
Elements of \lcm\ and \olcm\ are obtained from elements of \cm\ and
\lcm\ (respectively) by removing the upper boxes (thus converting
lower endpoints of upper curves to unmatched vertices).
For the case of \olcm, unmatched vertices will
still have an orientation (indicated by an arrow attached to the
vertex).  As for middle curves in the case of \ocm{\mathbf{a}_1, \dots, \mathbf{a}_l},
the unmatched
vertices in an element of \olcm\ must be arranged to that those
oriented down are to the right of those oriented up.
See Figure~\ref{fig:olcm}.

\begin{figure}
\begin{center}
\epsfig{file=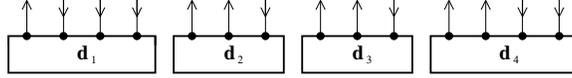,width=3in}
\caption{$\mathbf{d} = (4,3,3,4)$, $\mathbf{a} = (3,1,1,2)$ \label{fig:beta}}
\end{center}
\end{figure}

\begin{figure}
\begin{center}
\epsfig{file=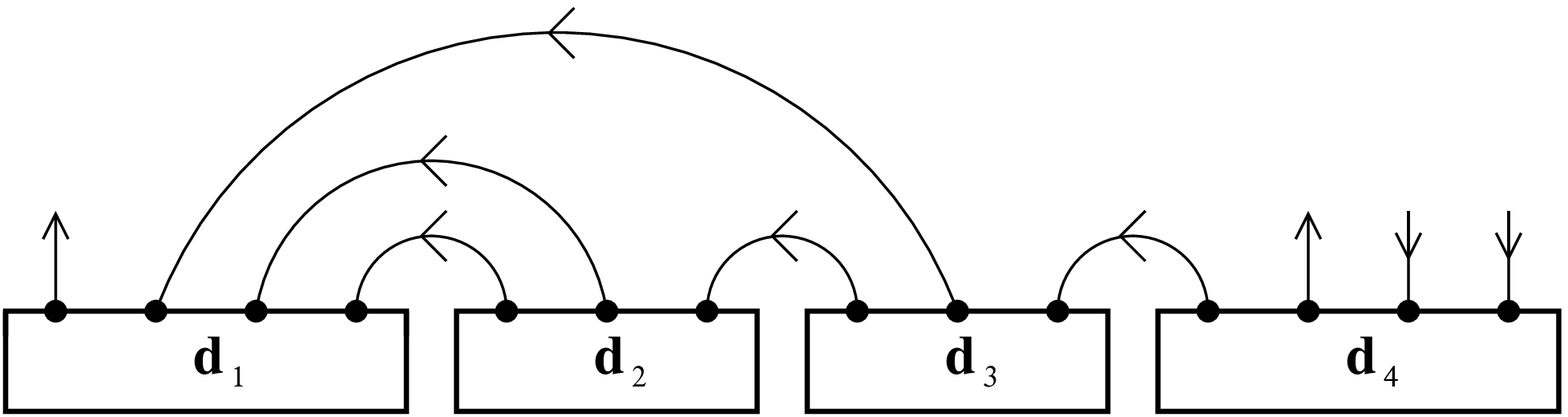,width=3in}
\caption{Oriented lower crossingless match associated to $\mathbf{d} =
    (4,3,3,4)$, $\mathbf{a} = (3,1,1,2)$ \label{fig:beta_cm}}
\end{center}
\end{figure}

Let $\mathbf{a} \in (\Z_{\ge 0})^k$ be such
that $\mathbf{a}_i \le \mathbf{d}_i$ for $i = 1,2,\ldots,k$.  We
associate an oriented lower crossingless match to $\mathbf{a}$ as follows.  For
each $i$, place down arrows on the rightmost $\mathbf{a}_i$ vertices of
the box representing $V_{\mathbf{d}_i}$.  Place up arrows on
the remaining vertices.  See Figure~\ref{fig:beta}.
There is a unique way to form an oriented lower crossingless match
such that the orientation of any curve agrees with the direction of
the arrows at its endpoints.  Namely,
starting from the right connect each down arrow to the first unmatched
up arrow to its right (if there is any).
Note that this produces an oriented lower crossingless match where the
unmatched vertices are arranged so that
all those with down arrows are to the right of those with up arrows
(otherwise, we could have matched more vertices).  See
Figure~\ref{fig:beta_cm}.
So to each $\mathbf{a}$ there is an associated element of \olcm.
Conversely, given an element of \olcm, there is exactly one
$\mathbf{a}$ which produces it.  So we have a one to one correspondence
between the set of elements $\mathbf{a} \in (\Z_{\ge 0})^k$
such that $\mathbf{a}_i \le \mathbf{d}_i$ and oriented lower
crossingless matches \olcm.  We will denote the oriented lower
crossingless match associated to $\mathbf{a}$ by $M(\mathbf{d}, \mathbf{a})$.

We can put a partial ordering on the sets \cm{\mathbf{a}_1, \dots, \mathbf{a}_l},
\ocm{\mathbf{a}_1, \dots, \mathbf{a}_l}, \lcm\ and \olcm\
as follows.  For any two elements $S_1$ and $S_2$ of one of these sets, $S_1
\le S_2$ if the set of lower curves of $S_1$ is a subset of the set of
lower curves of $S_2$.

Given the geometrization of irreducible representations of \U\
(Section \ref{sec:geom_reps}), it is
natural to seek a geometrization of the tensor product and the space
of intertwiners.  This geometric realization is the focus of Sections
\ref{sec:geom_tensor} and \ref{sec:intertwiners}.


\section{Geometric Realization of The Tensor Product}
\label{sec:geom_tensor}
\subsection{Definition of the Tensor Product Variety $\tv$}
\label{sec:tensor-variety}
We now describe a variety (introduced in
\cite{malkin} and \cite{nakajima:01})
corresponding to the tensor
product of the irreducible representations
$V_{\mathbf{d}_1},V_{\mathbf{d}_2},\ldots,V_{\mathbf{d}_k}$.
This construction will yield three distinct
bases of the tensor product in a natural way.

Fix a $d$-dimensional vector space $D$ and let $\mathbf{d} \in
(\Z_{\ge 0})^k$ be such that $\sum_{i=1}^k \mathbf{d}_i = d$.
Define
\begin{multline}
\tvo = \{(\mathbf{D}=\{\mathbf{D}_i\}_{i=0}^k, W)\ |\
0=\mathbf{D}_0\subset
\mathbf{D}_1\subset\cdots\subset \mathbf{D}_k=D,\, \\
W\subset D,\, \dim \mathbf{D}_i/\mathbf{D}_{i-1}=\mathbf{d}_i \}.
\end{multline}
\tvo\ admits a natural $GL(D)$ action.  Namely
\[
g \cdot (\{\mathbf{D}_i\}_{i=0}^k, W) = (\{g \mathbf{D}_i\}_{i=0}^k, gW)
\]
for $g \in GL(D)$ and $(\mathbf{D}, W) \in \tvo$.  Now let
\begin{multline}
\tv \stackrel{\text{def}}{=} \{(\mathbf{D}=\{\mathbf{D}_i\}_{i=0}^k, W, t)\ |\
0=\mathbf{D}_0\subset
\mathbf{D}_1\subset\cdots\subset \mathbf{D}_k=D,\, W\subset D,\, \\
t\in \End D,\, t(\mathbf{D}_i)\subset \mathbf{D}_{i-1},\, \dim \mathbf{D}_i/\mathbf{D}_{i-1}=\mathbf{d}_i,\, \im t
\subset W \subset \ker t\}.  
\end{multline}
We call \tv\ the \emph{tensor product variety}.
We say a flag
$\mathbf{D} = (0=\mathbf{D}_0\subset \mathbf{D}_1\subset \cdots \subset
\mathbf{D}_k=D)$ is \emph{$t$-stable} if
$t(\mathbf{D}_i)\subset \mathbf{D}_{i-1}$ for $i=1,\ldots,k$.

If we consider the corresponding varieties $\tvo'$ and $\tv'$ defined
over \barFq, a
straightforward computation shows that $\tv'$ is the union of the
conormal bundles of the orbits of the action of $GL(D)$ on $\tvo'$.

We define the action of $E$, $F$ and $K^{\pm 1}$ on the set of
functions on \tv\ just as for the other spaces considered so far.
Namely, let
\begin{equation*}
\tv[w;\mathbf{d}] = \{(\mathbf{D},W,t) \in \tv\, |\, \dim W = w\}
\end{equation*}
\begin{multline*}
\tv[w,w+1;\mathbf{d}] = \{(\mathbf{D},U,W,t)\, |\,
 (\mathbf{D},U,t),(\mathbf{D},W,t)\in \tv,\, U \subset W,\, \\
 \dim U=w,\, \dim W=w+1\}.
\end{multline*}
We then have the projections
\begin{equation}
\label{eq:tensor-EFaction}
\tv \stackrel{\pi_1}{\longleftarrow}
\bigcup_w \tv[w,w+1;\mathbf{d}]
\stackrel{\pi_2}{\longrightarrow} \tv.
\end{equation}
where $\pi_1(\mathbf{D},U,W,t) = (\mathbf{D},U,t)$ and $\pi_2(\mathbf{D},U,W,t) = (\mathbf{D},W,t)$.
The action of $E$, $F$ and $K^{\pm 1}$ is defined by \eqref{EF_action} as
usual.  Of course, the notation for the action of $K^{\pm 1}$ now means that
\begin{equation}
\label{eq:K_action_T}
(K^{\pm 1}f)(\mathbf{D},W,t) = q^{\pm(d - 2\dim W)} f(\mathbf{D},W,t).
\end{equation}


\subsection{A Set of Basic Functions on the Tensor Product Variety}
We now describe a set of basic functions on \tv\ which
will be used to form spaces of functions isomorphic to
$V_{\mathbf{d}_1}\otimes~\cdots~\otimes~V_{\mathbf{d}_k}$.
As usual, fix a $d$-dimensional vector space $D$.  For a flag
$\mathbf{D} = (0 = \mathbf{D}_0 \subset \dots \subset \mathbf{D}_k =
D)$ and a subspace $W
\subset D$, define
$\alpha(W,\mathbf{D}) \in (\Z_{\ge 0})^k$ by
\[
\alpha(W,\mathbf{D})_i = \dim (W \cap \mathbf{D}_i)/(W \cap \mathbf{D}_{i-1}).
\]
For
$\mathbf{w}, \mathbf{r}, \mathbf{n} \in(\Z_{\ge 0})^k$,
define
\begin{align}
\A = \{(\mathbf{D},W,t) \in \tv\, |\, & \alpha(W,\mathbf{D}) = \mathbf{w}, \\
& \alpha(\im t, \mathbf{D}) = \mathbf{r},\,
\alpha(\ker t, \mathbf{D}) = \mathbf{n}\}. \nonumber
\end{align}
Note that the non-empty sets $\A$ are precisely the orbits of the
action of $GL(D)$
given by
\[
g \cdot (\{\mathbf{D}_i\}_{i=0}^k, W, t) = (\{g
\mathbf{D}_i\}_{i=0}^k, gW, gtg^{-1}),\  g \in GL(D).
\]
From now on, the term \emph{constructible} will mean constructible
with respect to the stratification given by these sets.  We say that a
function $f$ on
\tv\ is \emph{invariant} if it is
invariant under the action of $GL(D)$ given by
\[
(g \cdot f)(x) = f(g^{-1} x),\; g\in GL(D).
\]
Let $\mathcal{T}(\mathbf{d})$ denote the space of invariant functions
on \tv.
We will also use the notation
\[
\mathbf{a}^{(j,l)} = \sum_{i=j}^l \mathbf{a}_i, \qquad
\s{\mathbf{a}} = \sum_{i=1}^k \mathbf{a}_i
\]
for $\mathbf{a} \in (\Z_{\ge 0})^k$
and we will let $\delta^j$ denote the element of $(\Z_{\ge 0})^k$ such
that $\delta^j_j
= 1$ and $\delta^j_i = 0$ for all $i \ne j$.

Let
\begin{equation}
\label{eq:k_def}
\K = q^{\sum_{i<j} (\mathbf{r}_i \mathbf{w}_j
  + \mathbf{w}_i \mathbf{n}_j - \mathbf{w}_i
  \mathbf{w}_j)}
\end{equation}
and define
\begin{equation}
\f = \K \F.
\end{equation}
Then it is easy to see that
\[
\mathcal{T}(\mathbf{d}) = \Span \{\f\}_{\mathbf{w}, \mathbf{r}, \mathbf{n}} .
\]
We will call the \f\ \emph{basic functions}.  Note that $\f = \F$ if
  $q=1$.  As will be seen below, the factor of
  $\K$ is neccessary in order for the $\f$ to correspond to certain
  vectors in the tensor product.   Note that $\f \equiv 0$ unless
  $\mathbf{r} \le \mathbf{w}
  \le \mathbf{n}$ where we define the partial ordering such that for
  $\mathbf{a}, \mathbf{b} \in (\Z_{\ge 0})^k$,
\begin{equation}
\label{eq:order_def}
\begin{split}
&\mathbf{a} \le \mathbf{b} \Longleftrightarrow \sum_{i=1}^j
\mathbf{a}_i \le \sum_{i=1}^j \mathbf{b}_i \mbox{ for } 1\le j
\le k. \\
&\mathbf{a} < \mathbf{b} \Longleftrightarrow \mathbf{a} \le
\mathbf{b},\, \mathbf{a} \ne \mathbf{b}
\end{split}
\end{equation}
Also, $\f \equiv 0$ unless $\s{\mathbf{r}} + \s{\mathbf{n}} =
\s{\mathbf{d}} = d$.

\begin{theorem}
\label{thm:basic_isom}
The action of \U\ described in
Section~\ref{sec:tensor-variety} endows $\mathcal{T}(\mathbf{d})$ with
the structure of a \U-module
and the map
\[
\eta_{\mathbf{r},\mathbf{n}} : \Span \{\f\}_{\mathbf{w}} \to
V_{\mathbf{n}_1 - \mathbf{r}_1} \otimes \dots \otimes V_{\mathbf{n}_k
  - \mathbf{r}_k}
\]
given by
\begin{equation}
\label{eq:basic_isom}
\eta_{\mathbf{r},\mathbf{n}}(\f) = \otimes \sv[\mathbf{n} - \mathbf{r}]{\mathbf{w} - \mathbf{r}}
\end{equation}
(and extended by linearity) is a \U -module isomorphism.
\end{theorem}

\begin{proof}
Fix a $(\mathbf{D},W,t) \in \tv$ such that $\alpha(W,\mathbf{D}) = \mathbf{w} - \delta^j$ for
some $j$ (it is easy to see that $E\f(\mathbf{D},W,t) = 0$ unless $W$
satisfies this property).  Then
\begin{align*}
E\f(\mathbf{D},W,t) &=  q^{-\dim (\pi_1^{-1}(\mathbf{D},W,t))}
(\pi_1)_!\pi_2^*\f(\mathbf{D},W,t)
\\
&= \K q^{-\dim (\pi_1^{-1}(\mathbf{D},W,t))} \chi_q
(\pi_1^{-1}(\mathbf{D},W,t) \cap
\pi_2^{-1}(\A)).
\end{align*}
Now,
\begin{align*}
\pi_1^{-1}(\mathbf{D},W,t) &\cong \{U\, |\, W \subset U \subset \ker t,\, \dim U =
\dim W + 1\} \\
&\cong \Proj^{\dim (\ker t) - \dim W - 1} \\
&= \Proj^{\s{\mathbf{n}} - (\s{\mathbf{w}}-1) - 1} \\
&= \Proj^{\s{\mathbf{n}} - \s{\mathbf{w}}}.
\end{align*}
So $\dim(\pi_1^{-1}(\mathbf{D},W,t)) = |\mathbf{n}| - |\mathbf{w}|$ and
\begin{align*}
\pi_1^{-1}(&\mathbf{D},W,t) \cap \pi_2^{-1}(\A) \\
&\cong \{U\, |\, W \subset U \subset \ker t,\, \alpha(U,\mathbf{D}) = \mathbf{w}\} \\
&\cong \{U\, |\, (W \cap \mathbf{D}_j)
\subset U \subset (\ker t \cap \mathbf{D}_j), \\
& \qquad \dim (U \cap \mathbf{D}_{j-1}) =
\mathbf{w}^{(1,j-1)},\,
\dim U = \mathbf{w}^{(1,j)}\} \\
&\cong \{U\, |\, U \subset (\ker t \cap \mathbf{D}_j)/(W \cap \mathbf{D}_j), \\
& \qquad U \not \subset
(\ker t \cap \mathbf{D}_{j-1})/(W \cap \mathbf{D}_{j-1}),\, \dim U = 1\} \\
&\cong \Proj^{\dim (\ker t \cap \mathbf{D}_j)/(W \cap \mathbf{D}_j) - 1} - \Proj^{\dim
  (\ker t \cap \mathbf{D}_{j-1})/(W \cap \ker \mathbf{D}_{j-1}) - 1} \\
&= \Proj^{\mathbf{n}^{(1,j)} - (\mathbf{w}-\delta^j)^{(1,j)} - 1} -
\Proj^{\mathbf{n}^{(1,j-1)} - (\mathbf{w} - \delta^j)^{(1,j-1)} - 1} \\
&= \Proj^{\mathbf{n}^{(1,j)} - \mathbf{w}^{(1,j)}} -
\Proj^{\mathbf{n}^{(1,j-1)} - \mathbf{w}^{(1,j-1)} - 1}.
\end{align*}
Thus
\begin{align*}
E\f(\mathbf{D},W,t) &= \K q^{-(\s{\mathbf{n}} - \s{\mathbf{w}})} \left(
  \sum_{i=0}^{\mathbf{n}^{(1,j)} - \mathbf{w}^{(1,j)}} q^{2i} -
  \sum_{i=0}^{\mathbf{n}^{(1,j-1)} - \mathbf{w}^{(1,j-1)} - 1} q^{2i}
  \right) \\
&= \K q^{\s{\mathbf{w}} - \s{\mathbf{n}}} \sum_{\mathbf{n}^{(1,j-1)} -
  \mathbf{w}^{(1,j-1)}}^{\mathbf{n}^{(1,j)} - \mathbf{w}^{(1,j)}} q^{2i} \\
&= \K q^{\s{\mathbf{w}} - \s{\mathbf{n}} + 2\left(\mathbf{n}^{(1,j-1)} -
  \mathbf{w}^{(1,j-1)}\right)} \sum_{i=0}^{\mathbf{n}_j -
  \mathbf{w}_j} q^{2i} \\
&= \K q^{-\mathbf{w}^{(1,j-1)} + \mathbf{w}^{(j+1,k)} + \mathbf{n}^{(1,j-1)} -
  \mathbf{n}^{(j+1,k)}} [\mathbf{n}_j - \mathbf{w}_j + 1].
\end{align*}
Now,
\[
\K[\mathbf{w} - \delta^j] = \K q^{-\mathbf{r}^{(1,j-1)} - \mathbf{n}^{(j+1,k)}
  + \mathbf{w}^{1,j-1} + \mathbf{w}^{j+1,k}}
\]
So
\begin{align*}
\K q^{-\mathbf{w}^{(1,j-1)} + \mathbf{w}^{(j+1,k)} + \mathbf{n}^{(1,j-1)} -
  \mathbf{n}^{(j+1,k)}} = \K[\mathbf{w}-\delta^j] q^{\mathbf{r}^{(1,j-1)} +
  \mathbf{n}^{(1,j-1)} - 2\mathbf{w}^{(1,j-1)}}
\end{align*}
and thus
\begin{align*}
E\f(\mathbf{D},W,t) = \K[\mathbf{w}-\delta^j] q^{\mathbf{r}^{(1,j-1)} +
  \mathbf{n}^{(1,j-1)} - 2\mathbf{w}^{(1,j-1)}} [\mathbf{n}_j -
  \mathbf{w}_j + 1].
\end{align*}
Therefore,
\begin{align}
E\f &= \sum_{j=1}^k  q^{\mathbf{r}^{(1,j-1)} +
  \mathbf{n}^{(1,j-1)} - 2\mathbf{w}^{(1,j-1)}} [\mathbf{n}_j -
  \mathbf{w}_j + 1]
  \K[\mathbf{w}- \delta^j]
  \id_{A_{\mathbf{w} - \delta^i, \mathbf{r}, \mathbf{n}}} \nonumber \\
&= \sum_{j=1}^k q^{\mathbf{r}^{(1,j-1)} +
  \mathbf{n}^{(1,j-1)} - 2\mathbf{w}^{(1,j-1)}} [\mathbf{n}_j -
  \mathbf{w}_j + 1]
  \f[\mathbf{w}- \delta^j] \nonumber \\
\label{eq:Ebasic}
&= \sum_{j=1}^k q^{\sum_{i=1}^{j-1} (\mathbf{n}_i - \mathbf{r}_i -
  2(\mathbf{w}_i - \mathbf{r}_i))} [\mathbf{n}_j - \mathbf{w}_j + 1]
  \f[\mathbf{w}- \delta^j].
\end{align}
Similarly
\begin{equation}
\label{eq:Fbasic}
F\f = \sum_{j=1}^k q^{-\sum_{i=j+1}^k (\mathbf{n}_i - \mathbf{r}_i -
  2(\mathbf{w}_i - \mathbf{r}_i))} [\mathbf{w}_j - \mathbf{r}_j
  +1]\f[\mathbf{w}+\delta^j].
\end{equation}
It follows immediately from \eqref{eq:K_action_T} that
\begin{equation}
\label{eq:Kbasic}
\begin{split}
K^{\pm 1} \f &= q^{\pm(d - 2\s{\mathbf{w}})} \f \\
&= q^{\pm \sum_{i=1}^k (\mathbf{n}_i - \mathbf{r}_i - 2(\mathbf{w}_i -
    \mathbf{r}_i))} \f
\end{split}
\end{equation}
since $\s{\mathbf{r}} + \s{\mathbf{n}} = \s{d}$.

Now recall that $x
\in \U$ acts on $V_{\mathbf{d}_1} \otimes \dots \otimes V_{\mathbf{d}_k}$ as
$\Delta^{(k-1)}(x)$.  In particular
\begin{equation}
\label{eq:EFK_coproduct}
\begin{split}
&\Delta^{(k-1)}E = \sum_{i=1}^k K\otimes \dots \otimes K \otimes E
\otimes 1 \otimes \dots \otimes 1 \\
&\Delta^{(k-1)}F = \sum_{i=1}^k 1 \otimes \dots \otimes 1 \otimes F
\otimes K^{-1} \otimes \dots \otimes K^{-1} \\
&\Delta^{(k-1)} K^{\pm 1} = K^{\pm 1} \otimes \dots \otimes K^{\pm 1}
\end{split}
\end{equation}
where in the first two equations, the $E$ or $F$ appears in the
$i^{th}$ position.
Comparing \eqref{eq:EFK_coproduct} and \eqref{U_action_V_d} to
\eqref{eq:Ebasic}, \eqref{eq:Fbasic} and, \eqref{eq:Kbasic} the result
follows.
\end{proof}


\subsection{The Space $\mathcal{T}_0(\mathbf{d})$ and the Elementary
  Basis $\B_e$}
\label{subsec:elem_basis}
Note that if $t=0$, then $\mathbf{r} = \mathbf{0}$ and $\mathbf{n} =
\mathbf{d}$.  Let $\B_e =
\{f_{\mathbf{w},\mathbf{0},\mathbf{d}}\}_\mathbf{w}$.  Then $\Span \B_e$ is the
space of invariant functions on $\tvo$ which we shall
denote by $\mathcal{T}_0(\mathbf{d})$.
We see from Theorem~\ref{thm:basic_isom} that the map
\begin{equation*}
\eta_{\mathbf{0},\mathbf{d}} : f_{\mathbf{w},\mathbf{0},\mathbf{d}}
\mapsto \otimes
\sv{\mathbf{w}}
\end{equation*}
(extended by linearity) is an isomorphism
\begin{equation*}
\mathcal{T}_0(\mathbf{d}) \cong  V_{\mathbf{d}_1}
\otimes \dots \otimes V_{\mathbf{d}_k}.
\end{equation*}
We have therefore exhibited the elementary basis as the set of
invariant functions
on the variety $\tvo \subset \tv$.


\subsection{The Space $\mathcal{T}_c(\mathbf{d})$}
\label{subsec:extensions}
The goal of this section is to develop a natural way to extend
invariant functions on
$\tvo$ to invariant functions on \tv\ with larger supports.  Recall that
$\mathcal{T}_0(\mathbf{d})$ and $\mathcal{T}(\mathbf{d})$ are the spaces
of invariant functions on
$\tvo$ and $\tv$ respectively.  The action of $E$,
$F$ and $K^{\pm 1}$ defined by \eqref{EF_action} gives both $\mathcal{T}_0(\mathbf{d})$ and
$\mathcal{T}(\mathbf{d})$ the structure
of a \U-module as can be seen from Theorem~\ref{thm:basic_isom}.
We will call a \U-module map $\epsilon : f \mapsto f^e$ from
$\mathcal{T}_0(\mathbf{d})$ to $\mathcal{T}(\mathbf{d})$ an
\emph{extension}.  Assuming an extension $\epsilon$ exists,
$\eta_{\mathbf{r},\mathbf{n}} \circ \epsilon \circ
(\eta_{\mathbf{0},\mathbf{d}})^{-1}$
is an intertwiner from $V_{\mathbf{d}_1} \otimes \dots \otimes
V_{\mathbf{d}_k}$ to $V_{\mathbf{n}_1 - \mathbf{r}_1} \otimes \dots \otimes
V_{\mathbf{n}_k - \mathbf{r}_k}$.  Conversely, each such set of
intertwiners determines an extension.  Namely, given a set
of intertwiners
\[
\{\gamma_{\mathbf{r}, \mathbf{n}}\ : V_{\mathbf{d}_1} \otimes
\dots \otimes V_{\mathbf{d}_k} \to V_{\mathbf{n}_1 - \mathbf{r}_1} \otimes \dots
\otimes V_{\mathbf{n}_k - \mathbf{r}_k}\}_{\mathbf{r}, \mathbf{n}}
\]
we extend a function $f \in \mathcal{T}_0(\mathbf{d})$ to a function
$f^e \in \mathcal{T}(\mathbf{d})$
by defining
\begin{equation}
\label{eq:extension}
f^e = \sum_{\mathbf{r}, \mathbf{n}} (\eta_{\mathbf{r}, \mathbf{n}})^{-1} \circ
\gamma_{\mathbf{r}, \mathbf{n}} \circ \eta_{\mathbf{0}, \mathbf{d}} \; (f) .
\end{equation}

\begin{figure}
\begin{center}
\epsfig{file=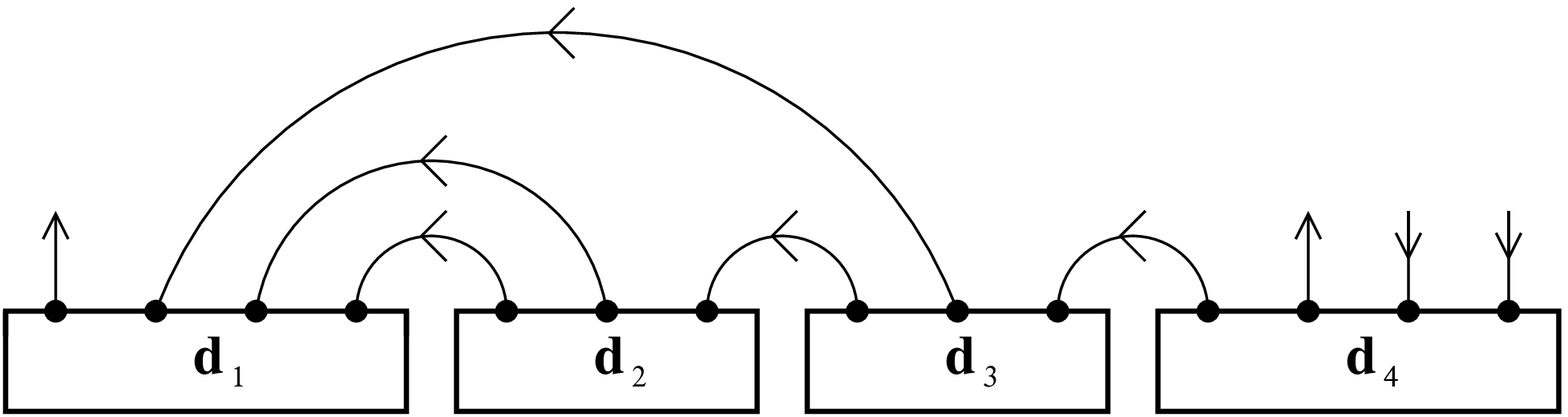,width=3in}
\caption{Oriented lower crossingless match S.  $\mathbf{r}^S = (3,1,1,0)$, $\mathbf{n}^S =
(4,1,1,3)$ \label{fig:ext1}}
\end{center}
\end{figure}

\begin{figure}
\begin{center}
\epsfig{file=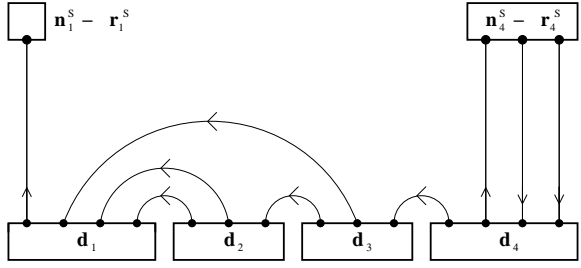,width=3in}
\caption{Completion of Figure~\ref{fig:ext1} to an oriented
  crossingless match to
$V_{\mathbf{n}^S_1 - \mathbf{r}^S_1} \otimes \dots \otimes
V_{\mathbf{n}^S_k - \mathbf{r}^S_k} = V_1 \otimes V_0 \otimes V_0 \otimes
V_3$ \label{fig:ext2}}
\end{center}
\end{figure}

From Section~\ref{sec:crossmatch} we know that a basis for the space of
intertwiners between two tensor product representations of \U\ is
given by the corresponding crossingless matches.  Now, a lower curve
represents a particular action of $t \in \End W$.  A lower curve connecting
$V_{\mathbf{d}_i}$ and $V_{\mathbf{d}_j}$ with $i < j$ represents the fact
that $t$ sends a vector in $\mathbf{D}_j - \mathbf{D}_{j-1}$ to a vector in $\mathbf{D}_i -
\mathbf{D}_{i-1}$.
So for any lower
crossingless match $S$, fix a basis of $D$ compatible with the flag
$\mathbf{D}$ and let $t$ be the map whose matrix in this basis has
$(i,j)$ component equal to 1 if $i<j$ and $S$ has an curve connecting
the $i^{th}$ and $j^{th}$ vertices and is equal to zero otherwise.
Then let $\mathbf{r}^S$ and $\mathbf{n}^S$ be defined
as $\alpha(\im t, \mathbf{D})$ and $\alpha(\ker t, \mathbf{D})$.
Thus, $\mathbf{r}^S_i$ is the number of left endpoints
of the lower curves contained in $V_{\mathbf{d}_i}$ and $\mathbf{n}^S_i$ is
$\mathbf{d}_i$ minus the number of right endpoints of the lower curves
contained in $V_{\mathbf{d}_i}$.  See Figure~\ref{fig:ext1}.
Then complete $S$ to a crossingless
match to $V_{\mathbf{n}^S_1 - \mathbf{r}^S_1} \otimes \dots \otimes
V_{\mathbf{n}^S_k - \mathbf{r}^S_k}$  as in Figure~\ref{fig:ext2} (there is a
unique way to do this).
Let
${\tilde \gamma}_{\mathbf{r}^S, \mathbf{n}^S}$ be the corresponding
intertwiner in the dual basis (that is, commuting with the action
of \U\ given by ${\bar \Delta}^{(k-1)}$).
Note that
${\tilde \gamma}_{\mathbf{r}^S, \mathbf{n}^S}$ is well-defined since the map
$S \mapsto (\mathbf{r}^S, \mathbf{n}^S)$ described above is injective.
Now let $\gamma_{\mathbf{r}^S, \mathbf{n}^S} = \sigma {\tilde
  \gamma}_{\mathbf{r}^S,
  \mathbf{n}^S} \sigma$.  As noted in
Section~\ref{sec:crossmatch},
$\gamma_{\mathbf{r}^S, \mathbf{n}^S} : V_{\mathbf{d}_1} \otimes
\dots \otimes
V_{\mathbf{d}_k} \to V_{\mathbf{n}^S_1 - \mathbf{r}^S_1} \otimes \dots \otimes
V_{\mathbf{n}^S_k - \mathbf{r}^S_k}$ is an intertwiner in the usual
basis (that is, it commutes with the action of \U\ given by
$\Delta^{(k-1)}$).  For all $(\mathbf{r},
\mathbf{n})$ not of the form $(\mathbf{r}^S, \mathbf{n}^S)$ for some lower
crossingless match $S$, let
$\gamma_{\mathbf{r}, \mathbf{n}} = 0$.  Then let $\epsilon : f \mapsto
f^e$ be the map defined by \eqref{eq:extension}.

\begin{prop}
\label{prop:epsilon}
The extension $\epsilon$ is an isomorphism onto its image and
\[
f^e |_{\tvo} = f.
\]
\end{prop}

\begin{proof}
This follows immediately from
Theorem~\ref{thm:basic_isom}.
\end{proof}

Let
\[
\mathcal{T}_c(\mathbf{d}) = \epsilon (\mathcal{T}_0(\mathbf{d}))
\subset \mathcal{T}(\mathbf{d}) .
\]
It follows from Proposition~\ref{prop:epsilon} and
Theorem~\ref{thm:basic_isom} that
$\mathcal{T}_c(\mathbf{d}) \cong V_{\mathbf{d}_1} \otimes \dots 
\otimes V_{\mathbf{d}_k}$.  And it follows from
Proposition~\ref{prop:epsilon} that $\epsilon :
\mathcal{T}_0(\mathbf{d}) \to \mathcal{T}_c(\mathbf{d})$ is an
isomorphism of \U-modules with inverse given by restriction to \tvo.
We will find a distinguished basis of $\mathcal{T}_c(\mathbf{d})$
related to the irreducible components of $\tv'$.  Before we do this,
we must first examine these irreducible components.


\subsection{The Irreducible Components of the Tensor Product Variety}
For the
remainder of this section we consider varieties defined over \barFq.
To avoid confusion,
we denote the corresponding varieties by $\tvo'$ and $\tv'$.
For $\mathbf{w} \in (\Z_{\ge 0})^k$ such that
$\mathbf{w}_i \le \mathbf{d}_i$, let
\begin{equation*}
Z_\mathbf{w}' = \{(\mathbf{D},W,t) \in \tv'\, |\, \alpha(W,\mathbf{D})
= \mathbf{w}\}.
\end{equation*}
We then have the following.

\begin{theorem}
$\{\overline{Z_\mathbf{w}'}\}_\mathbf{w}$ are the irreducible
components of $\tv'$.
\end{theorem}

\begin{proof}
It is obvious that $\sqcup_\mathbf{w} Z_\mathbf{w}' = \tv'$ (where
$\sqcup$ denotes
disjoint union).  Also, the connected components of $\tv'$ are given by
fixing the dimension of $W$. Thus, since $\s{\mathbf{w}} = \dim W$, it
suffices to prove that the $Z_\mathbf{w}'$ are
irreducible and locally closed and that $\dim Z_\mathbf{w}'$ is independent
of $\mathbf{w}$ for fixed $\s{\mathbf{w}}$.  Consider the maps
\begin{equation*}
Z_\mathbf{w}' \stackrel{p_1}{\rightarrow} {^1}Z_{\mathbf{w}}'
\stackrel{p_2}{\rightarrow} {^2}Z_{\mathbf{w}}'
\end{equation*}
where
\begin{align*}
{^1}Z_{\mathbf{w}}' &= \{(\mathbf{D},W)\, |\, (\mathbf{D},W,t) \in Z_\mathbf{w}'
\mbox{ for some } t\} \\
{^2}Z_{\mathbf{w}}' &= \{\mathbf{D}\, |\, (\mathbf{D},W) \in {^1}Z_\mathbf{w}' \mbox{
  for some } W\} \\
p_1(\mathbf{D},W,t) &= (\mathbf{D},W) \\
p_2(\mathbf{D},W) &= \mathbf{D}.
\end{align*}
Then $p_1$ and $p_2$ are locally trivial fibrations.  Now
\begin{equation*}
{^2}Z_{\mathbf{w}}' = \{ \mathbf{D}=\{\mathbf{D}_i\}_{i=0}^k\, |\, 0 = \mathbf{D}_0 \subset \mathbf{D}_1 \subset
\dots \subset \mathbf{D}_k=D,\, \dim \mathbf{D}_i/\mathbf{D}_{i-1} = \mathbf{d}_i\}
\end{equation*}
is simply a flag manifold.  It is a homogeneous space as follows.
$GL(D)$ acts transitively on ${^2}Z_{\mathbf{w}}'$ with stabilizer isomorphic to
the set of matrices
\begin{equation*}
G_0 = \left\{ \left.
\begin{pmatrix}
M_1    & \ast   & \cdots & \ast    \\
0      & M_2    & \ddots & \vdots  \\
\vdots & \ddots & \ddots & \ast    \\
0      & \cdots & 0      & M_k
\end{pmatrix}
\right|\, M_i \in GL(\mathbf{d}_i) \right\}.
\end{equation*}
Thus,
\begin{equation}
\begin{split}
\label{eq:basedim}
\dim {^2}Z_{\mathbf{w}}' &= \dim GL(D) - \dim G_0 \\
&= \sum_{i < j} \mathbf{d}_i \mathbf{d}_j
\end{split}
\end{equation}
Now, the fiber of $p_2$ over a point $\mathbf{D} \in {^2}Z_{\mathbf{w}}'$ is
\begin{equation*}
F_2= \{ W \subset D\, |\, \alpha(W,\mathbf{D}) = \mathbf{w} \}.
\end{equation*}
The group $G_0$ acts transitively on this space and the stabilizer is
isomorphic to the set of matrices
\begin{equation*}
G_1 = \left\{ \left.
\begin{pmatrix}
M_1     & \ast   & \ast   & \ast   & \ast   & \cdots & \cdots & \ast  \\
0       & N_1    & 0      & \ast   & 0      & \cdots & \cdots & \ast  \\
0       & 0      & M_2    & \ast   & \ast   &        &        & \ast  \\
0       & 0      & 0      & N_2    & 0      &        &        & \ast  \\
0       & 0      & 0      & 0      & \ddots & \ddots &        & \ast  \\
0       & 0      & 0      & 0      & \ddots & \ddots & 0      & \vdots
\\
\vdots  & \vdots & \vdots & \vdots & \ddots & \ddots & M_k    & \ast  \\
0       & 0      & 0      & 0      & \cdots & \cdots & 0      & N_k
\end{pmatrix}
\right|,
\begin{matrix} M_i \in GL(\mathbf{w}_i), \\
               N_i \in GL(\mathbf{d}_i - \mathbf{w}_i)
\end{matrix}
\right\}.
\end{equation*}
So
\begin{equation}
\begin{split}
\label{eq:F2dim}
\dim F_2 &= \dim G_0 - \dim G_1 \\
&= \sum_{j \le i} \mathbf{w}_i(\mathbf{d}_j - \mathbf{w}_j)
\end{split}
\end{equation}

The fiber of $p_1$ over a point $(\mathbf{D},W) \in {^1}Z_{\mathbf{w}}'$ is
\begin{equation}
\label{eq:F1_def}
F_1 = \{ t \in \End D\, |\, t(\mathbf{D}_i) \subset
\mathbf{D}_{i-1},\, \im t \subset W \subset \ker t\}.
\end{equation}
Pick a basis $\{ u_i \}_{i=1}^d$ of $D$ such that $\{ u_i
\}_{i=1}^{\mathbf{d}_1 + \dots + \mathbf{d}_j}$ is a basis for
$\mathbf{D}_j$ and
\[
\bigcup_{l=0}^j \{ u_i \}_{i= \mathbf{d}^{(1,l-1)}+
  1}^{\mathbf{d}^{(1,l-1)} +
  \mathbf{w}_l}
\]
(where $\mathbf{d}^{(1,0)}=0$) is a basis for $W \cap
\mathbf{D}_j$.  Then by considering the
matrices of $t$ in this basis it is easy to see that $F_1$ is an
affine space of dimension
\begin{equation}
\label{eq:F1dim}
\dim F_1 = \sum_{i<j} \mathbf{w}_i(\mathbf{d}_j -
\mathbf{w}_j).
\end{equation}
So from Equations \eqref{eq:basedim}, \eqref{eq:F2dim} and
\eqref{eq:F1dim} we see that
\begin{equation}
\begin{split}
\dim Z_\mathbf{w}' &= \sum_{i < j} \mathbf{d}_i \mathbf{d}_j + \sum_{i,j=1}^k
\mathbf{w}_i(\mathbf{d}_j - \mathbf{w}_j) \\
&= \sum_{i < j} \mathbf{d}_i \mathbf{d}_j + \s{\mathbf{w}}(d - \s{\mathbf{w}})
\end{split}
\end{equation}
and thus $\dim Z_\mathbf{w}'$ is independent of $\mathbf{w}$ (for a fixed value
of $\s{\mathbf{w}}$).

Now, ${^2}Z_{\mathbf{w}}'$, $F_1$ and $F_2$ are all smooth and connected,
hence irreducible.  Also, ${^2}Z_{\mathbf{w}}'$ and $F_1$ are closed while
$F_2$ is locally closed.  The latter statement follows from the fact
that $F_2$ is equal to the closed set $\{W \subset D\, |\,
\alpha(W, \mathbf{D}) \ge \mathbf{w}\}$ minus the finite collection of
closed sets $\{W \subset D\, |\, \alpha(W, \mathbf{D}) \ge
\mathbf{a}\}_{\mathbf{a} > \mathbf{w}}$.
Thus each $Z_\mathbf{w}'$ is irreducible and locally closed.
\end{proof}

Let
\begin{align}
\A' = \{(\mathbf{D},W,t) \in \tv'\, |\, & \alpha(W,\mathbf{D}) = \mathbf{w}, \\
& \alpha(\im t, \mathbf{D}) = \mathbf{r},\,
\alpha(\ker t, \mathbf{D}) = \mathbf{n}\}. \nonumber
\end{align}
We will need the following two propositions in the sequel.

\begin{prop}
\label{prop:A_irred_relation_1}
Let $\mathbf{w} \in (\Z_{\ge 0})^k$ with $\mathbf{w}_i \le
\mathbf{d}_i$ for all $i$ and let $M=M(\mathbf{d},\mathbf{w})$.
Then $A_{\mathbf{w},
  \mathbf{r}^M, \mathbf{n}^M}'$ is an open dense
subset of $\overline{Z_\mathbf{w}'}$.  In particular, $\overline
{A_{\mathbf{w}, \mathbf{r}^M, \mathbf{n}^M}'} = \overline{Z_\mathbf{w}'}$.
\end{prop}

\begin{proof}
It is enough to show that $A_{\mathbf{w}, \mathbf{r}^M,
  \mathbf{n}^M}'$ is dense 
in $Z_\mathbf{w}'$ (it is obvious from the definitions that $A_{\mathbf{w},
  \mathbf{r}^M, \mathbf{n}^M}' \subset Z_\mathbf{w}'$).
Since $\mathbf{r}^M \le \mathbf{w} \le \mathbf{n}^M$ by construction of $M$, we
have that the projection of $A_{\mathbf{w}, \mathbf{r}^M,
    \mathbf{n}^M}'$ onto ${^1}Z_\mathbf{w}'$ is all of ${^1}Z_\mathbf{w}'$.  Thus it
suffices to show that $A_{\mathbf{w}, \mathbf{r}^M, \mathbf{n}^M}'$
is dense in each fiber.  Fix $(\mathbf{D},W) \in {^1}Z_\mathbf{w}'$.
The fiber, $F_1$, of the projection $p_1$ is given by \eqref{eq:F1_def}.  The
intersection of $F_1$ with $(p_1 |_{A_{\mathbf{w}, \mathbf{r}^M,
    \mathbf{n}^M}'})^{-1} (\mathbf{D},W)$ is isomorphic to
\begin{multline*}
B = \{t \in \End D\, |\, t(\mathbf{D}_i) \subset \mathbf{D}_{i-1},\,
\im t \subset W \subset \ker t, \\
\alpha(\ker t, \mathbf{D}) = \mathbf{n}^M,\,
\alpha(\im t, \mathbf{D}) = \mathbf{r}^M \}.
\end{multline*}
Choose a basis $\beta$ of $D$ compatible with the flag $\mathbf{D}$ and
subspace $W$ (that is, there exist bases for $W$ and each $\mathbf{D}_i$
which are subsets of $\beta$).
Now, since $\im t \subset W \subset \ker t$, $t$ can be factored
through $D/W$ and considered as a map into $W$.  Each $t$ is uniquely
determined by the corresponding ${\bar t} \in \End(D/W, W)$.  Consider
the matrix of ${\bar t}$ in the basis of $D/W$ given by the projection
of the basis $\beta$ under the natural map $D \to D/W$ and the basis
of $W$ which is a subset of $\beta$.  It must be of the following form:
\begin{equation}
C_t = 
\begin{pmatrix}
0      & A_{1,2} & A_{1,3} & \cdots & A_{1,k}   \\
\vdots & 0       & A_{2,3} & \cdots & A_{2,k}   \\
\vdots & \vdots  & \ddots  & \ddots & \vdots    \\
\vdots & \vdots  & \ddots  & 0      & A_{k-1,k} \\
0      & 0       & 0       & 0      & 0         \\
\end{pmatrix}
\end{equation}
where $A_{i,j}$ is a $(\mathbf{w}_i) \times (\mathbf{d}_j -
\mathbf{w}_j)$ matrix.
Then $t \in B$ if and only if each submatrix
\begin{equation}
\label{eq:submatrix}
C_t^{i,j} = 
\begin{pmatrix}
A_{i,i+1} & A_{i,i+2}   & \cdots & A_{i,j+1}   \\
0         & A_{i+1,i+2} & \cdots & A_{i+1,j+1} \\
\vdots    & \ddots      & \ddots & \vdots    \\
0         & \cdots      & 0      & A_{j,j+1} \\
\end{pmatrix}
\ , \quad 1 \le i \le j \le k-1 ,
\end{equation}
has maximal rank.
To see this, consider the diagram $M'$ of non-crossing oriented curves
connecting the
$V_{\mathbf{d}_i}$ associated to a $t \in F_1$.  That is, the number
of down-oriented vertices among those associated to $V_{\mathbf{d}_i}$
is given by $\mathbf{w}_i$ and the number of
left and right
endpoints of curves of $M'$ in $V_{\mathbf{d}_i}$ are given by
$\alpha(\im t, \mathbf{D})_i$ and $\mathbf{d}_i - \alpha(\ker t, \mathbf{D})_i$
respectively.  A priori, this is not an oriented lower crossingless
match (for instance, the unmatched vertices of $M'$ might not be
arranged so that those oriented down are to the right of those oriented
up).
The requirement that $C_t^{i,j}$ has
maximal rank is equivalent to the requirement that $M'$ has the
maximum possible number of curves connecting $V_{\mathbf{d}_{i}}$,
$V_{\mathbf{d}_{i+1}}$, \dots, and $V_{\mathbf{d}_{j+1}}$.  Thus,
referring to the
definition of $M(\mathbf{d}, \mathbf{w})$ given in
Section~\ref{sec:crossmatch}, we see that the condition that all the
$C_t^{i,j}$ have 
maximal rank is equivalent to the condition that $M' = M$ (so
$M'$ is indeed an oriented crossingless match) and thus equivalent to
$\alpha(\im t, \mathbf{D}) = \mathbf{r}^{M'} = \mathbf{r}^M$ and
$\alpha(\ker t, \mathbf{D}) = \mathbf{n}^{M'} = \mathbf{n}^M$ or $t
\in B$.  Note that this argument also allows us to see that $B$ is not
empty since it contains the element $t$ given by the matrix whose
$(i,j)$ entry is 1
if $i<j$ and $M$ contains a curve connecting the $i^{th}$
and $j^{th}$ vertices and zero otherwise.  In fact, this is a
canonical form of any $t \in B$.  That is, by a change of basis
(preserving the flag $\mathbf{D}$), we
can transform the matrix of any $t \in B$ to this form.

Assume we know
that the subset $N_{m,n}$ of the set
$M_{m,n}$ of $m \times n$ matrices given by
\[
N_{m,n} = \{A \in M_{m,n} \, |\, A \text{ has maximal rank}\} 
\]
is an open subset of the set $M_{m,n}$.  Then $N_{m,n}$ is given by the
non-vanishing of a finite collection of polynomials in the matrix elements of
$M_{m,n}$ (recall we are working in the Zariski topology).  Thus, the
requirement that the submatrices $C_t^{i,j}$ have maximal
rank is equivalent to the non-vanishing of a finite number of
polynomials in the matrix elements of the $C_t^{i,j}$ (and hence of
$C_t$).  Therefore,
we will have shown that $B$
is the intersection of a finite number of open subsets of $F_1$ and
hence is open (and thus dense since it is not empty) in $F_1$.

So it remains to show that $N_{m,n}$ is dense in $M_{m,n}$.  But if we
let $r=\min(m,n)$, then
\[
N_{m,n} = \{A \in M_{m,n}\, |\, \text{At least one $r \times r$
  submatrix of $A$ has rank $r$}\}.
\]
which is a union of open subsets of $M_{m,n}$ (since an $r \times r$
matrix has rank $r$ if and only if its determinant is non-zero) and
hence open (and dense) in $M_{m,n}$.
\end{proof}

\begin{prop}
\label{prop:A_irred_relation_2}
With the notation of Proposition \ref{prop:A_irred_relation_1}, $\overline
{A_{\mathbf{a}, \mathbf{r}^S, \mathbf{n}^S}'} \subset \overline{Z_\mathbf{w}'}$
for all $S \le M$, $\mathbf{a} \ge \mathbf{w}$, $\s{\mathbf{a}} = \s{\mathbf{w}}$.
\end{prop}

\begin{proof}
It suffices to show that $A_{\mathbf{a}, \mathbf{r}^S, \mathbf{n}^S}'$
is contained
in $\overline{Z_\mathbf{w}'}$.  The image of $A_{\mathbf{a}, \mathbf{r}^S,
  \mathbf{n}^S}'$ under the projection $p_1$ is $\{(\mathbf{D},W)\, |\,
\alpha(W,\mathbf{D}) = \mathbf{a} \}$ which is contained in
$\overline{{^1}Z_\mathbf{w}'}$ since $\mathbf{a} \ge \mathbf{w}$ and
$\s{\mathbf{a}} = \s{\mathbf{w}}$.  The
fiber of the projection $p_1$ (restricted to $A_{\mathbf{a},
  \mathbf{r}^S, \mathbf{n}^S}'$) over a point $(\mathbf{D},W)$ is
\[
\{t\, |\, (\mathbf{D},W,t) \in \tv,\, \alpha(\ker t, \mathbf{D}) =
\mathbf{n}^S,\, \alpha(\im t,\mathbf{D}) = \mathbf{r}^S \}
\]
and this is in the closure of the set
\[
\{t\, |\, (\mathbf{D},W,t) \in \tv,\, \alpha(\ker t, \mathbf{D}) =
\mathbf{n}^M,\, \alpha(\im t,\mathbf{D}) = \mathbf{r}^M \}
\]
since $S \le M$.  So $A_{\mathbf{a}, \mathbf{r}^S, \mathbf{n}^S}' \subset
\overline{Z_\mathbf{w}'}$.
\end{proof}

We now define the irreducible components of \tv\ to be the \Fq\ points
of the irreducible components of $\tv'$.  Let
$\overline{Z_\mathbf{w}}$ denote the set of \Fq\ points of the irreducible
component $\overline{Z_\mathbf{w}'}$ of $\tv'$.  We also define the \emph{dense
  points} of an irreducible component $\overline{Z_\mathbf{w}}$ of
\tv\ to be the \Fq\ points of
the dense subset $A_{\mathbf{w}, \mathbf{r}^M, \mathbf{n}^M}'$ (where $M =
M(\mathbf{d}, \mathbf{w})$) of the corresponding irreducible component
$\overline{Z_\mathbf{w}'}$ of $\tv'$.  However, the \Fq\ points of
$\A'$ are exactly the elements of $\A$.  Thus, the dense points of the
irreducible component $\overline{Z_\mathbf{w}}$ of
\tv\ are just the points of $A_{\mathbf{w}, \mathbf{r}^M, \mathbf{n}^M}$.


\subsection{Geometric Realization of the Canonical Basis}
\label{subsec:can_basis}

We are now ready to describe the set of functions mentioned at the end
of Section~\ref{subsec:extensions}.
Define
\begin{equation}
\label{def:h_and_g}
\begin{split}
&h^\mathbf{d}_\mathbf{w} = \eta_{\mathbf{0},\mathbf{d}}^{-1} \left(
  \diamondsuit \sv{\mathbf{w}} \right) \\
&g^\mathbf{d}_\mathbf{w} = \left( h^\mathbf{d}_\mathbf{w} \right)^e \\
&\B_c = \{ g^\mathbf{d}_\mathbf{w} \}_\mathbf{w}.
\end{split}
\end{equation}

For a vector $w_1 \otimes \dots \otimes w_k \in V_{\mathbf{a}_1}
\otimes \dots V_{\mathbf{a}_k}$, let $(w_1 \otimes \dots \otimes
w_k)^r = w_k \otimes \dots \otimes w_1 \in V_{\mathbf{a}_k}
\otimes \dots V_{\mathbf{a}_1}$.
For an intertwiner $\gamma : V_{\mathbf{a}_1} \otimes \dots \otimes
V_{\mathbf{a}_k} \to V_{\mathbf{b}_1} \otimes \dots \otimes V_{\mathbf{b}_l}$
corresponding to a
crossingless match $S$, let $\gamma^\dagger : V_{\mathbf{b}_l} \otimes
\dots \otimes V_{\mathbf{b}_1} \to V_{\mathbf{a}_k} \otimes \dots \otimes
V_{\mathbf{a}_1}$ denote the
intertwiner corresponding to the crossingless match $S$ rotated $180^\circ$.  
It follows easily from the graphical calculus described in
\cite{frenkel-khovanov:97} that
\[
\pair{\gamma(v)}{w} = \pair{v}{(\sigma \gamma^\dagger \sigma) (w)} =
\pair{ v}{(\tilde\gamma)^\dagger (w)}
\]
for any $v \in V_{\mathbf{a}_1} \otimes \dots \otimes V_{\mathbf{a}_k}$ and $w \in
V_{\mathbf{b}_l} \otimes \dots \otimes V_{\mathbf{b}_1}$.

We will need the following results.

\begin{lemma}
\label{lemma:can_to_can}
\[
\gamma_{\mathbf{r}^S, \mathbf{n}^S} \left(\diamondsuit \sv{\mathbf{w}}
  \right) =
\begin{cases}
\diamondsuit \sv[\mathbf{n}^S - \mathbf{r}^S]{\mathbf{w} - \mathbf{r}^S}
  & \text{if $S \le M(\mathbf{d},\mathbf{w})$}, \\
0 & \text{otherwise}.
\end{cases}
\]
\end{lemma}
\begin{proof}
It is apparent from the graphical calculus of
\cite{frenkel-khovanov:97} that if $S \le M(\mathbf{d},\mathbf{w})$,
then $\left( {\tilde \gamma}_{\mathbf{r}^S,
  \mathbf{n}^S} \right)^\dagger \left( \left( \heartsuit \dv[\mathbf{n}^S -
  \mathbf{r}^S]{\mathbf{w} - \mathbf{r}^S}
  \right)^r \right) = \left( \heartsuit \dv{\mathbf{w}}
\right)^r$ and that $\left( {\tilde \gamma}_{\mathbf{r}^S, \mathbf{n}^S}
\right)^\dagger$ sends
other dual canonical basis elements $\left( \heartsuit \dv[\mathbf{n}^S -
    \mathbf{r}^S]{\mathbf{a}} \right)^r$, $\mathbf{a} \ne \mathbf{w} - \mathbf{r}^S$, to
elements of the form $\left( \heartsuit \dv{\mathbf{a}'}
\right)^r$ with $\mathbf{a}' \ne \mathbf{w}$. Therefore
\begin{align*}
\pair{\gamma_{\mathbf{r}^S, \mathbf{n}^S} \left( \diamondsuit
    \sv{\mathbf{w}} \right)}{\left( \heartsuit \dv[\mathbf{n}^S -
    \mathbf{r}^S]{\mathbf{w} - \mathbf{r}^S} \right)^r} &= \pair{\diamondsuit
    \sv{\mathbf{w}}}{\left( {\tilde \gamma}_{\mathbf{r}^S,
    \mathbf{n}^S} \right)^\dagger \left( \left( \heartsuit \dv[\mathbf{n}^S -
    \mathbf{r}^S]{\mathbf{w} - \mathbf{r}^S} \right)^r \right)}\\
&= \pair{\diamondsuit
    \sv{\mathbf{w}}}{\left( \heartsuit
    \dv{\mathbf{w}} \right)^r} \\
&= 1
\end{align*}
and
\[
\pair{\gamma_{\mathbf{r}^S, \mathbf{n}^S} \left( \diamondsuit
    \sv{\mathbf{w}} \right)}{\left( \heartsuit \dv[\mathbf{n}^S -
    \mathbf{r}^S]{\mathbf{a}} \right)^r} = 0
\]
for all $\mathbf{a} \ne \mathbf{w} - \mathbf{r}^S$.  Thus $\gamma_{\mathbf{r}^S,
  \mathbf{n}^S} \left(\diamondsuit \sv{\mathbf{w}} \right) =
\diamondsuit \sv[\mathbf{n}^S - \mathbf{r}^S]{\mathbf{w} - \mathbf{r}^S}$.  A similar
  argument demonstrates that $\gamma_{\mathbf{r}^S, \mathbf{n}^S}
  \left(\diamondsuit \sv{\mathbf{w}} \right) = 0$ if $S \not \le
  M(\mathbf{d}, \mathbf{w})$ since then the image of $\left( {\tilde
  \gamma}_{\mathbf{r}^S, \mathbf{n}^S} \right)^\dagger$ is spanned by
  $\heartsuit \dv{\mathbf{a}}$ with $\mathbf{a} \ne \mathbf{w}$.
\end{proof}

\begin{prop}
\label{prop:g_form}
\[
g^{\mathbf{d}}_\mathbf{w} = \sum_{S \le M(\mathbf{d}, \mathbf{w})} \left(
  \eta_{\mathbf{r}^S \mathbf{n}^S} \right)^{-1} \left( \diamondsuit \sv[\mathbf{n}^S -
  \mathbf{r}^S]{\mathbf{w} - \mathbf{r}^S} \right).
\]
\end{prop}

\begin{proof}
This follows immediately from Lemma~\ref{lemma:can_to_can}.
\end{proof}

\begin{prop}
\label{prop:can_basis}
$\diamondsuit \sv{\mathbf{w}}$ is equal to $\otimes
  \sv{\mathbf{w}}$ plus a linear combination of
  elements $\otimes \sv{\mathbf{a}}$, $\mathbf{a} >
  \mathbf{w}$, $\s{\mathbf{a}} = \s{\mathbf{w}}$, with coefficients
  in $q^{-1}\N[q^{-1}]$.
\end{prop}

\begin{proof}
This follows from Sections 1.5 and 1.6 of \cite{frenkel-khovanov:97}.
\end{proof}

We can now prove one of our main results.

\begin{theorem}
\label{thm:B_c}
$g^\mathbf{d}_\mathbf{w}$ is the unique element of
$\mathcal{T}_c(\mathbf{d})$, up to a multiplicative constant,
satisfying the following conditions.
\begin{enumerate}
\item \label{cond1} $g^{\mathbf{d}}_\mathbf{w}$ is equal to a non-zero constant
on the set of dense
points $A_{\mathbf{w}, \mathbf{r}^M, \mathbf{n}^M}$ of the irreducible
component $\overline{Z_\mathbf{w}}$ (where $M=M(\mathbf{d}, \mathbf{w})$).
\item \label{cond2} The support of $g^{\mathbf{d}}_\mathbf{w}$ lies in
  $\overline{Z_\mathbf{w}}$.
\end{enumerate}
Furthermore, the set $\{g^\mathbf{d}_\mathbf{w}\}_\mathbf{w}$ is a
basis of $\mathcal{T}_c(\mathbf{d})$ and the map
\[
\diamondsuit \sv{\mathbf{w}} \mapsto g^\mathbf{d}_\mathbf{w}
\]
(extended by linearity) is a \U -module isomorphism $V_{\mathbf{d}_1}
\otimes \dots \otimes V_{\mathbf{d}_k} \cong \mathcal{T}_c(\mathbf{d})$.
\end{theorem}

\begin{proof}
In this proof, to simplify notation in calculations, we will
suppress the
isomorphism $\eta_{\mathbf{r}, \mathbf{n}}$ defined by
\eqref{eq:basic_isom} and identify the vector $\otimes \sv[\mathbf{n} -
  \mathbf{r}]{\mathbf{w} -\mathbf{r}}$ with the function $\f$.
Assume that ${g'}^\mathbf{d}_\mathbf{w}$ satisfies the above
conditions and let ${g'}^\mathbf{d}_\mathbf{w} = \left(
  {h'}^\mathbf{d}_\mathbf{w} \right)^e$.
The value of
${g'}^{\mathbf{d}}_\mathbf{w}$ on $A_{\mathbf{w}, \mathbf{r}^M, \mathbf{n}^M}$ is
given by $k_{\mathbf{w}, \mathbf{r}^M, \mathbf{n}^M}$ times the coefficient of
$f_{\mathbf{w}, \mathbf{r}^M, \mathbf{n}^M}$ when
${g'}^{\mathbf{d}}_\mathbf{w}$ is written as a linear combination of the basic
functions.  This coefficient is equal to
\begin{equation*}
\pair{\gamma_{\mathbf{r}^M, \mathbf{n}^M}
  ({h'}^{\mathbf{d}}_\mathbf{w})}{\left( \otimes \dv[\mathbf{n}^M -
  \mathbf{r}^M]{\mathbf{w} - \mathbf{r}^M} \right)^r}
\end{equation*}.

Therefore, since $k_{\mathbf{w}, \mathbf{r}^M, \mathbf{n}^M} \ne 0$,
condition~\ref{cond1} is equivalent to 
\begin{align*}
\pair{\gamma_{\mathbf{r}^M, \mathbf{n}^M}
  ({h'}^{\mathbf{d}}_\mathbf{w})}{\left( \otimes \dv[\mathbf{n}^M -
  \mathbf{r}^M]{\mathbf{w} - \mathbf{r}^M} \right)^r} \ne 0 \\
\Leftrightarrow \pair{
  {h'}^{\mathbf{d}}_\mathbf{w}}{({\tilde \gamma}_{\mathbf{r}^M, \mathbf{n}^M})^\dagger
  \left( \left( \otimes \dv[\mathbf{n}^M -
  \mathbf{r}^M]{\mathbf{w} - \mathbf{r}^M} \right)^r \right)} \ne 0.
\end{align*}

Now, since $M=M(\mathbf{d}, \mathbf{w})$ is the oriented lower
crossingless match
associated to $\mathbf{w}$, $M(\mathbf{n}^M -
\mathbf{r}^M, \mathbf{w} - \mathbf{r}^M)$ has no lower curves and
all down arrows are to the right of all up arrows.  So after being rotated
by $180^\circ$ (but keeping the original orientation of unmatched
vertices -- for example, those oriented up remain oriented up), this
diagram has all down arrows to the left of all up arrows.
Thus, by Section~2.3 of \cite{frenkel-khovanov:97}, $\left( \otimes
\dv[\mathbf{n}^M - \mathbf{r}^M]{\mathbf{w} - \mathbf{r}^M}
\right)^r = \left( \heartsuit
\dv[\mathbf{n}^M - \mathbf{r}^M]{\mathbf{w} - \mathbf{r}^M}
\right)^r$.  It also follows from the graphical calculus of
\cite{frenkel-khovanov:97} that
\[
({\tilde \gamma}_{\mathbf{r}^M, \mathbf{n}^M})^\dagger
  \left( \left( \heartsuit \dv[\mathbf{n}^M -
  \mathbf{r}^M]{\mathbf{w} - \mathbf{r}^M} \right)^r \right) = 
  \left( \heartsuit \dv{\mathbf{w}} \right)^r.
\]
Therefore condition~\ref{cond1} is equivalent to
\begin{equation}
\label{eq:g_proof_1}
\pair{
  {h'}^{\mathbf{d}}_\mathbf{w}}{\left( \heartsuit \dv{\mathbf{w}}
  \right)^r} \ne 0.
\end{equation}

Next we consider condition~\ref{cond2}.  In order for this condition
to be satisfied, ${g'}^{\mathbf{d}}_{\mathbf{w}}$ must be equal to zero
on $A_{\mathbf{w}',
  \mathbf{r}^{M'}, \mathbf{n}^{M'}}$ for all $\mathbf{w}' \ne \mathbf{w}$ (where $M'
  = M(\mathbf{d}, \mathbf{w}')$).  By an argument analogous to
that given above, this is equivalent to the condition
\begin{equation}
\label{eq:g_proof_2}
\pair{
  {h'}^{\mathbf{d}}_\mathbf{w}}{\left( \heartsuit \dv{\mathbf{w}'}
  \right)^r} = 0
\end{equation}
for all $\mathbf{w}' \ne \mathbf{w}$.  Therefore, by
\eqref{eq:g_proof_1} and \eqref{eq:g_proof_2}, we must have
\begin{equation*}
{h'}^{\mathbf{d}}_\mathbf{w} = c^{\mathbf{d}}_\mathbf{w} \cdot
\diamondsuit \sv{\mathbf{w}} = c^{\mathbf{d}}_\mathbf{w} \cdot
h^{\mathbf{d}}_\mathbf{w}
\end{equation*}
for some non-zero constant $c^{\mathbf{d}}_\mathbf{w}$ which proves
uniqueness up to a multiplicative constant.
It still remains to show that $g^{\mathbf{d}}_\mathbf{w}$ satisfies
the given conditions.

Now, by Proposition~\ref{prop:g_form}, the value of
$g^{\mathbf{d}}_\mathbf{w}$ on $A_{\mathbf{w}, \mathbf{r}^M,
  \mathbf{n}^M}$, where $M =
M(\mathbf{d}, \mathbf{w})$, is equal to \K\ times the coeffient of $\otimes
  \sv[\mathbf{n}^M - \mathbf{r}^M]{\mathbf{w}-\mathbf{r}^M}$ in the expression of
\linebreak
$\diamondsuit
  \sv[\mathbf{n}^M - \mathbf{r}^M]{\mathbf{w}-\mathbf{r}^M}$ as a linear
combination of elementary
basis elements.  By Proposition~\ref{prop:can_basis}, this coefficient
is equal to
1.  So $g^{\mathbf{d}}_\mathbf{w}$ is equal to a non-zero constant on
$A_{\mathbf{w} ,\mathbf{r}^M, \mathbf{n}^M}$.
Also, by Propositions \ref{prop:g_form} and \ref{prop:can_basis},
$g^{\mathbf{d}}_\mathbf{w}$ is equal to a linear combination of
functions of the form $\left( \eta_{\mathbf{r}^S, \mathbf{n}^S} \right)^{-1}
\left( \otimes \sv[\mathbf{n}^S - \mathbf{r}^S]{\mathbf{a}} \right) =
f_{\mathbf{a} +
\mathbf{r}^S, \mathbf{r}^S, \mathbf{n}^S}$ with $S \le M$,
$\s{\mathbf{a}} = \s{\mathbf{w} -
  \mathbf{r}^S}$ ($\Rightarrow \s{\mathbf{a} + \mathbf{r}^S} =
\s{\mathbf{w}}$), and $\mathbf{a} \ge \mathbf{w} - \mathbf{r}^S$
($\Rightarrow \mathbf{a} + \mathbf{r}^s \ge \mathbf{w}$).  Thus, by
Proposition~\ref{prop:A_irred_relation_2}, the support of
$g^{\mathbf{d}}_\mathbf{w}$ lies in $\overline{Z_\mathbf{w}}$.
So we have demonstrated that the functions
$g^{\mathbf{d}}_{\mathbf{w}}$ are
the unique functions, up to a
multiplicative constant, satisfying conditions \ref{cond1} and
\ref{cond2}.

The last two statements of the theorem follow from the fact that the
map $\eta_{\mathbf{0},\mathbf{d}} : \diamondsuit
\sv{\mathbf{w}} \mapsto
h^{\mathbf{d}}_{\mathbf{w}}$ (extended by linearity) is
a \U -module isomorphism 
$V_{\mathbf{d}_1} \otimes \dots \otimes V_{\mathbf{d}_k} \cong
\mathcal{T}_0(\mathbf{d})$  and the fact that $\epsilon$ is an isomorphism onto
its image.
\end{proof}


\subsection{A Conjectured Characterization of
  $\mathcal{T}_c(\mathbf{d})$ and $\B_c$}
\label{sec:perv-conj}

We present here a conjecture concerning an alternative characterization of
the basis $\B_c$.  Let \Perv\ be the category of semisimple perverse
sheaves on $\tv'$
constructible with respect to the stratification given by the $\A'$ and
let $\mathcal{D} : \Perv \to \Perv$ be the operation of Verdier
Duality.  For $\mathbf{B}^{\bullet} \in \Perv$ and $x \in \tv'$,
$\mathbf{B}_x^{\bullet}$
denotes the stalk complex at the point $x$.
We define the action of the involution $\Psi^{(k)}$ on
$\mathcal{T}(\mathbf{d})$ by
\[
\Psi^{(k)} \left( \f \right) = (\eta_{\mathbf{r}, \mathbf{n}})^{-1} \Psi^{(k)}
\eta_{\mathbf{r}, \mathbf{n}} \left( \f \right) ,
\]
where on the right hand side $\Psi^{(k)}$ is the involution used to
characterize the canonical basis (see \cite{frenkel-khovanov:97},
Section~1.6).  In
particular, the canonical basis is invariant under the action of
$\Psi^{(k)}$.  Now let $\theta : \Perv \to \mathcal{T}(\mathbf{d})$ be the map such that
\[
(\theta(\mathbf{B}^{\bullet}))(x) = \sum_i (-1)^i q^i (\mathbf{B}_x^i)
\text{ for $x \in \tv$ and $\mathbf{B}^{\bullet} \in \Perv$}.
\]

For each irreducible component $\overline{Z_\mathbf{w}'}$ of $\tv'$ there is
an intersection sheaf complex $IC_\mathbf{w}^{\bullet}$ associated to the
local system which is
the constant sheaf $\C$ (in degree zero) on the dense subset $A_{\mathbf{w},
  \mathbf{r}^M,
  \mathbf{n}^M}'$ where $M = M(\mathbf{d}, \mathbf{w})$ (see
\cite{kirwan} for details).

\begin{conjecture}
\label{conj:perv}
$\theta \left( IC_{\mathbf{w}}^{\bullet} \right) =
k_{\mathbf{w}, \mathbf{r}^M,
  \mathbf{n}^M}^{-1} g^{\mathbf{d}}_{\mathbf{w}}$.
\end{conjecture}
The factor of $k_{\mathbf{w}, \mathbf{r}^M,
  \mathbf{n}^M}^{-1}$ arises from the fact
that $\theta \left( IC_{\mathbf{w}}^{\bullet} \right)$ is equal to
one on the set $A_{\mathbf{w}, \mathbf{r}^M,
  \mathbf{n}^M}'$.
The proof of Conjecture~\ref{conj:perv} would most likely center around the
idea that the action of Verdier Duality in \Perv\ should correspond to
the action of $\Psi^{(k)}$ in $\mathcal{T}(\mathbf{d})$.  The precise
  statement is the following:
\begin{conjecture}
$\theta \D = \Psi^{(k)} \theta$.
\end{conjecture}


\section{Geometric Realization of The Intertwiners}
\label{sec:intertwiners}
\subsection{Defining the Intertwiners}
The goal of this section is to decompose \tv\ into subsets
corresponding to a basis for the space of intertwiners
\begin{equation}
\twine{\mu} = \Hom_{\U} (V_{\mathbf{d}_1} \otimes \cdots \otimes
V_{\mathbf{d}_k}, V_\mu).
\end{equation}
Note that $\twine{\mu} = 0$ unless $\mu = d-2r$ for some $0 \le r \le
d/2$ (where $d = \s{\mathbf{d}}$).
Thus, the intertwiners
will be maps from $\mathcal{T}(\mathbf{d})$ to
$\mathcal{M}(d)$ since these $V_\mu$ are precisely the
representations appearing in $\mathcal{M}(d)$ (see
Section~\ref{sec:geom_reps}).

Let $Y$ be a constructible subset of $\tv$. Define $R_Y :
\mathcal{T}(\mathbf{d}) \to \mathcal{T}(\mathbf{d})$
to be the map
which restricts
functions to their values on $Y$.  That is, for
$f \in \mathcal{T}(\mathbf{d})$, $R_Yf=\id_Y f$ (where the multiplication of
functions is pointwise).

Consider the map
$p : \tv \to \qrept[d]$
such that
$p(\mathbf{D},W,t)=(W,t)$.
Let $T_Y = p_! R_Y$.  Then $T_Y$ is a
map from $\mathcal{T}(\mathbf{d})$ to $\mathcal{M}(d)$.

\begin{prop}
\label{prop:TY}
If $Y \subset \tv$ satisfies $\pi_1 \pi_2^{-1} (Y) \subset Y$ and $\pi_2
\pi_1^{-1} (Y) \subset Y$ where $\pi_1$ and $\pi_2$ are the maps from
\eqref{eq:tensor-EFaction} then $T_Y$ is an intertwiner.
\end{prop}

\begin{proof}
It suffices to show that $T_Y$ commutes with the action of $E$,
$F$ and $K^{\pm 1}$ since these elements generate \U.
Note that the condition $\pi_1 \pi_2^{-1} (Y) \subset Y$ implies $\pi_2^{-1}(Y)
\subset \pi_1^{-1}(Y)$ and the condition $\pi_2 \pi_1^{-1} (Y) \subset
Y$ implies $\pi_1^{-1}(Y) \subset \pi_2^{-1}(Y)$.  Thus $\pi_1^{-1}(Y)
= \pi_2^{-1} (Y)$.
We first show that $T_Y E = E T_Y$.  Now $T_Y E = p_! R_Y E$ and
$E T_Y = E p_! R_Y$.  Thus it suffices to show that $R_Y E = E
R_Y$ and $p_! E = E p_!$.  Since 
$\mathcal{T}(\mathbf{d})$ is spanned by functions of the form $\id_A$
where $A$ is a
subvariety of \tv, we
need only check that actions agree on such functions.  For $x =
(\mathbf{D},W,t) \in \tv$
\begin{align*}
R_Y E \id_A (x) &= \id_Y(x) (E\id_A)(x) \\
&= \id_Y(x) q^{-\dim (\pi_1^{-1} (x))} ((\pi_1)_! \pi_2^* \id_A)(x) \\
&= q^{-\dim (\pi_1^{-1} (x))} \id_Y(x) ((\pi_1)_! \id_{\pi_2^{-1}(A)})(x) \\
&= q^{-\dim (\pi_1^{-1} (x))} \id_Y(x) \chi_q(\pi_1^{-1}(x) \cap
\pi_2^{-1}(A)) \\
&= q^{-\dim (\pi_1^{-1} (x))} \chi_q(\pi_1^{-1}(x \cap Y) \cap
\pi_2^{-1}(A)) \\
&= q^{-\dim (\pi_1^{-1} (x))} \chi_q(\pi_1^{-1}(x) \cap
\pi_1^{-1}(Y) \cap \pi_2^{-1}(A)) \\
&= q^{-\dim (\pi_1^{-1} (x))} \chi_q(\pi_1^{-1}(x) \cap
\pi_2^{-1}(Y) \cap \pi_2^{-1}(A)) \\
&= q^{-\dim (\pi_1^{-1} (x))} \chi_q(\pi_1^{-1}(x) \cap
\pi_2^{-1}(Y \cap A)) \\
&= q^{-\dim (\pi_1^{-1} (x))} (\pi_1)_!\id_{\pi_2^{-1}(Y \cap A)}(x) \\
&= q^{-\dim (\pi_1^{-1} (x))} (\pi_1)_!\pi_2^*\id_{Y \cap A}(x) \\
&= q^{-\dim (\pi_1^{-1} (x))} (\pi_1)_!\pi_2^*(\id_Y \id_A)(x) \\
&= E R_Y \id_A (x),
\end{align*}
where the fifth equality holds from consideration of the two cases $x
\in Y$ and $x \not \in Y$.

It remains to show that $p_!E = Ep_!$.  For the purposes of
this demonstration, we introduce the map
\[
p' : \bigcup_w \tv[w,w+1;\mathbf{d}] \rightarrow
\bigcup_w \dqrept{w}{w+1}{d}
\]
which acts as
$p'(\mathbf{D},U,W,t) = (U,W,t)$.  We have the following commutative
diagram.
\begin{equation*}
\begin{CD}
\tv @>{p}>> \qrept[d] \\
         @AA{\pi_2}A                                 @AA{\pi_2}A \\
\bigcup_w \tv[w,w+1;\mathbf{d}] @>{p'}>>
         \bigcup_w \dqrept{w}{w+1}{d} \\
         @VV{\pi_1}V                                 @VV{\pi_1}V \\
\tv @>{p}>> \qrept[d]
\end{CD}
\end{equation*}
As before we use the notation $\pi_1$ and $\pi_2$ to denote
several different, but analogous maps.

Note that $p^{-1} \pi_2 = \pi_2 (p')^{-1}$ (both are the map
$(U,W,t) \mapsto \{(\mathbf{D},W',t')\in \tv\, |\, W'=W, t'=t\}$).  Using this fact we
show that $(p')_!\pi_2^* = \pi_2^*p_!$.  Let $x \in
\tv$.  Then
\begin{align*}
(\pi_2^* p_! \id_A)(x) &= (p_!\id_A)(\pi_2(x)) \\
&= \chi_q(p^{-1}(\pi_2(x)) \cap A) \\
&= \chi_q(\pi_2((p')^{-1}(x)) \cap A) \\
&= \chi_q(\pi_2((p')^{-1}(x) \cap \pi_2^{-1}(A))) \\
&= \chi_q((p')^{-1}(x) \cap \pi_2^{-1}(A)) \\
&= (p')_! \id_{\pi_2^{-1}(A)}(x) \\
&= (p')_! \pi_2^* \id_A (x)
\end{align*}
where in the fourth equality we used the general fact that $\pi_2(B) \cap A =
\pi_2 (B \cap \pi_2^{-1}(A))$ and in the fifth equality we used the
fact that if $x=(U',W',t')$ then
\begin{align*}
(p')^{-1}(x) \cap \pi_2^{-1}(A) &= \{(\mathbf{D},U,W,t)\, |\,
p'(\mathbf{D},U,W,t)=(U',W',t'), \\
& \qquad \qquad \qquad \qquad \qquad \qquad \qquad
\pi_2(\mathbf{D},U,W,t) \in A\} \\
&= \{(\mathbf{D},U,W,t)\, |\, U=U',\, W=W',\, t=t',\, (\mathbf{D},W,t) \in A\} \\
&\cong \{(\mathbf{D},W,t)\, |\, W=W',\, t=t',\, (\mathbf{D},W,t) \in A\} \\
&= \pi_2((p')^{-1}(x) \cap \pi_2^{-1}(A)).
\end{align*}
We also have that $\pi_1 p' = p \pi_1$ (both are the map
$(\mathbf{D},U,W,t) \mapsto (U,t)$).  Thus
\begin{align*}
Ep_! &= q^{-\dim (\pi_1^{-1} (\cdot))} (\pi_1)_! \pi_2^* p_! \\
&= q^{-\dim (\pi_1^{-1} (\cdot))} (\pi_1)_! (p')_! \pi_2^* \\
&= q^{-\dim (\pi_1^{-1} (\cdot))} (\pi_1 p')_! \pi_2^* \\
&= q^{-\dim (\pi_1^{-1} (\cdot))} (p \pi_1)_! \pi_2^* \\
&= q^{-\dim (\pi_1^{-1} (\cdot))} p_!(\pi_1)_! \pi_2^* \\
&= p_! q^{-\dim (\pi_1^{-1} (\cdot))} (\pi_1)_! \pi_2^* \\
&= p_! E
\end{align*}
where we have used the fact that the map $f \mapsto f_!$ is functorial
\cite{macpherson:74}.

Thus, we have shown that $T_Y E = E T_Y$.  The proof that $T_Y F = F
T_Y$ is analogous.  Also,
\begin{align*}
K^{\pm 1} T_Y f(\mathbf{D},W,t) &= q^{\pm (d-2\dim W)} T_Y f(\mathbf{D},W,t) \\
&= T_Y q^{\pm (d-2\dim W)} f(\mathbf{D},W,t) \\
&= T_Y K^{\pm 1} f(\mathbf{D},W,t).
\end{align*}
\end{proof}


\subsection{A Basis $\B_I$ for the Space of Intertwiners}
\label{sec:intertwiner_basis}
We see from Section~\ref{sec:crossmatch} that a basis for the space of
intertwiners \twine{\mu}\ is in one to one correspondence with the set
of crossingless matches $\cm{\mu}$.
Note that crossingless matches of the form $\cm{\mu}$ (i.e. with only
one box on the top vertical line) are in one to one correspondence
with elements of \lcm.  For a given element S of \lcm\ simply set $\mu$
equal to the number of unmatched vertices of S and join the unmatched
vertices to the upper box.  Recall that elements $\mathbf{a} \in
(\Z_{\ge 0})^n$ such that $\mathbf{a}_i~\le~\mathbf{d}_i$ are in one to
one correspondence with the elements of \olcm.
Given such an $\mathbf{a}$, consider its associated oriented lower crossingless
match $M(\mathbf{d},\mathbf{a})$.  Note that $\s{\mathbf{a}}$
is the number of vertices (both matched and unmatched) in
$M(\mathbf{d},\mathbf{a})$ which are
oriented down.

For any flag $\mathbf{D}$ and $t \in \End D$ let
$\alpha(t,\mathbf{D}) = \alpha(\ker t,\mathbf{D})$.
Then let
\begin{equation}
\label{eq:Y_beta_def}
\begin{split}
Y_\mathbf{a} &= \{(\mathbf{D},W,t) \in \tv\, |\, \alpha(t,\mathbf{D}) =
\mathbf{n}^{M(\mathbf{d},\, \mathbf{a})},\, \dim W = \s{\mathbf{a}}\} \\
&= \bigcup_{\mathbf{w} : \s{\mathbf{w}} = \s{\mathbf{a}}} \bigcup_{\mathbf{r}} A_{\mathbf{w},
  \mathbf{r}, \mathbf{n}^{M(\mathbf{d}, \mathbf{a})}}.
\end{split}
\end{equation}

Now, note that $\mathbf{n}^{M(\mathbf{d}, \mathbf{a})}$ depends only on the
lower curves of $\mathbf{a}$ and not on the orientation of the unmatched
vertices.  Thus, if $\bar{\mathbf{a}}$
denotes the (unoriented) lower crossingless match associated to
$\mathbf{a}$, we can unambiguously define $\mathbf{n}^{\bar{\mathbf{a}}} =
\mathbf{n}^{M(\mathbf{d}, \mathbf{a})}$.
Then if $b$ is an unoriented crossingless match, we define
\begin{equation}
\label{eq:Y_b_def}
\begin{split}
Y_b &= \{(\mathbf{D},W,t) \in \tv\, |\, \alpha(t,\mathbf{D}) = \mathbf{n}^b \} \\
&= \bigcup_{\mathbf{a} : {\bar{\mathbf{a}}} = b} Y_\mathbf{a}.
\end{split}
\end{equation}
The last equality arises from the fact that if $(\mathbf{D},W,t) \in
\tv$ then $\im t \subset W \subset \ker t$, so $r \le \dim W \le
d - r$ (where $r= \rank t$).  Thus, since $(\mathbf{D},W,t) \in Y_b$
implies that $r = \rank t$ is the number of
lower curves in $b$, the values
$r,r+1, \ldots, d - r$ are precisely the number of down arrows
(that is, the $\s{\mathbf{a}}$) in the various $\mathbf{a}$ such that
$\bar{\mathbf{a}} = b$.  We also have the
following:

\begin{prop}
\label{prop:Y_b_union}
$\sqcup_b Y_b = \sqcup_\mathbf{a} Y_\mathbf{a} = \tv$ .
\end{prop}

\begin{proof}
It is obvious that the $Y_\mathbf{a}$ are disjoint.  Thus, from equation
\eqref{eq:Y_b_def} we see that it suffices to prove
that for every $(\mathbf{D},W,t) \in \tv$, $\alpha(t,\mathbf{D}) =
\mathbf{n}^{M(\mathbf{d}, \mathbf{a})}$ for some
crossingless match $\mathbf{a}$.  Fix an $(\mathbf{D},W,t) \in \tv$ and let $\mathbf{a} =
\alpha(t,\mathbf{D})$.  Now, down arrows of $\mathbf{a}$ represent
dimensions of the
kernel of $t$ while up arrows of $\mathbf{a}$ represent dimensions of
$D/ \ker t$.  Let $c$ denote the $i^{th}$ up arrow from the
left.  Since $\im t \subset \ker t$ and $t(\mathbf{D}_j) \subset(\mathbf{D}_{j-1})$,
there must be at least $i$ down arrows to the left of
$c$.  Since this holds for all $i$, it follows that each up arrow of
$M(\mathbf{d},\mathbf{a})$ is matched.  Thus, since
$\mathbf{n}^{M(\mathbf{d}, \mathbf{a})}$ is obtained from
$\mathbf{a}$ by forcing all unmatched vertices to be oriented down, we have
that $\mathbf{n}^{M(\mathbf{d}, \mathbf{a})} = \mathbf{a} =
\alpha(t,\mathbf{D})$.
\end{proof}

Define
\begin{equation}
\B_I = \left\{ T_{Y_b}\ \left|\, b \in \bigcup_{\mu} \cm{\mu} \right. \right\}.
\end{equation}

\begin{prop}
\label{prop:B_I}
Each element of $\B_I$ is an intertwiner and
\[
T_{Y_b}(\mathcal{T}(\mathbf{d}))
\subset \mathcal{M}^r(d) \cong V_\mu
\]
for $b \in \cm{\mu}$ and $r=(d-\mu)/2$.
\end{prop}

\begin{proof}
According to
Proposition~\ref{prop:TY}, to show that $T_{Y_b}$ is an intertwiner we
need only check that $\pi_2 \pi_1^{-1}
(Y_b) \subset
Y_b$ and $\pi_1 \pi_2^{-1} (Y_b) \subset Y_b$ for all
$b \in \cm{\mu}$.  If we
denote by $t^x$ and $\mathbf{D}^x$ the map $t$ and flag $\mathbf{D}$
of the point $x\in
\tv$ (so $x=(\mathbf{D}^x,W,t^x)$ for some $W$), then $t^y
= t^x$ and $\mathbf{D}^y = \mathbf{D}^x$ for all $y \in \pi_2 \pi_1^{-1} (x)$.  Thus
$\alpha(t^x,\mathbf{D}^x) =
\alpha(t^y,\mathbf{D}^y)$ for all $y \in \pi_2 \pi_1^{-1} (x)$ which
implies that $\pi_2 \pi_1^{-1}
(Y_b) \subset Y_b$ for all $b$.  Similarly $\pi_1 \pi_2^{-1} (Y_b)
\subset Y_b$ for all $b$.
Now, the image of $T_{Y_b}$ consists of functions on
\qreptr{d}\ where $r$ is the number of lower curves in $b$.
In fact, it is easy to see that for $f \in \mathcal{T}(\mathbf{d})$,
$T_{Y_b}(f)(W,t)$ depends only on the dimension of $W$ and the rank of
$t$.  So the image of $T_{Y_b}$ is contained in $\mathcal{M}^r(d)$.
Recall from
Section~\ref{sec:geom_reps}
that $\mathcal{M}^r(d) \cong V_{d-2r}$.  Since $r$ is equal
to the number of lower curves in $b$,
$d-2r$ is equal to the number of middle curves and hence
$d-2r=\mu$. So $T_{Y_b}$ is an
intertwiner into the representation $V_\mu$ as it should be.
\end{proof}


\subsection{The Space $\mathcal{T}_s(\mathbf{d})$ and the Basis $\B_s$}
\label{sec:B_d_def}
For the purposes of this section we will identify the sets \lcm\ and
$\cup_\mu \cm{\mu}$ as in Section~\ref{sec:intertwiner_basis}.  Also,
to simplify notation, we shall identify elements $\mathbf{a} \in
(\Z_{\ge 0})^k$ such that $\mathbf{a}_i \le \mathbf{d}_i$ with their
associated oriented lower crossingless matches
$M(\mathbf{d},\mathbf{a})$.

Let $\mathcal{T}_s(\mathbf{d})$ be the space of all functions $f \in
\mathcal{T}(\mathbf{d})$ such that
\[
\dim W = \dim W',\ \alpha(t, \mathbf{D}) = \alpha(t',
\mathbf{D}') \Rightarrow
f(\mathbf{D}, W, t) = f(\mathbf{D}', W', t').
\]
It is obvious that if we define
\begin{equation}
\B_s = \left\{ \id_{Y_\mathbf{a}}\, \left| \, \mathbf{a} \in
\bigcup_\mu \ocm{\mu} \right. \right\} ,
\end{equation}
then
\[
\mathcal{T}_s(\mathbf{d}) = \Span \B_s .
\]

\begin{theorem}
\label{thm:nat-basis}
$\mathcal{T}_s(\mathbf{d})$ is isomorphic as a \U -module to
$V_{\mathbf{d}_1} \otimes \cdots \otimes V_{\mathbf{d}_k}$ and $\B_s$
is a basis for $\mathcal{T}_s(\mathbf{d})$ adapted
to its decomposition into a direct sum of irreducible
representations.  That is, for a given $b \in \cm{\mu}$, the space
$\Span \{\id_{Y_\mathbf{a}}\, |\, \bar{\mathbf{a}} = b\}$ is isomorphic to the irreducible
representation $V_\mu$ via the map
\[
\id_{Y_\mathbf{a}} \mapsto \sv[\mu]{\mu - 2(\text{\# of unmatched down
    arrows in $\mathbf{a}$})}
\]
(extended by linearity).
\end{theorem}

\begin{proof}
For $\mathbf{a} \in \lcm$ such that $\mathbf{a}$ has at least one
unmatched up arrow, let $\mathbf{a}^+$ be the element of \lcm\
obtained from $\mathbf{a}$ by switching the orientation of the
rightmost unmatched up arrow.  Thus $\overline {\mathbf{a}^+} =
\bar{\mathbf{a}}$ and $\mathbf{a}^+$ has one more unmatched down arrow than
$\mathbf{a}$.  Similary, if $\mathbf{a} \in \lcm$ has at least one
unmatched down arrow, let $\mathbf{a}^-$ be the element of \lcm\
obtained from $\mathbf{a}$ by switching the orientation of the
leftmost unmatched down arrow.
Recall from the proofs of Propositions~\ref{prop:TY} and
\ref{prop:B_I} that
$\pi_1^{-1}(Y_b) = \pi_2^{-1}(Y_b)$.  It follows from this and the
fact that $Y_b = \bigcup_{\mathbf{a} : \bar{\mathbf{a}} = b}
Y_\mathbf{a}$ that $\pi_2
\pi_1^{-1} (Y_\mathbf{a}) = Y_{\mathbf{a}^+}$ if $\mathbf{a}$ has at
least one unmatched up arrow and $\pi_2
\pi_1^{-1} (Y_\mathbf{a}) = \emptyset$ otherwise.  Similarly, $\pi_1
\pi_2^{-1} (Y_\mathbf{a}) = Y_{\mathbf{a}^-}$ if $\mathbf{a}$ has at
least one unmatched down arrow and
$\pi_1 \pi_2^{-1} (Y_\mathbf{a}) = \emptyset$ otherwise.

Now, for $x \in \tv$,
\begin{align*}
F\id_{Y_\mathbf{a}}(x) &= q^{-\dim (\pi_2^{-1}(x))} (\pi_2)_! \pi_1^*
\id_{Y_\mathbf{a}} (x) \\
&= q^{-\dim (\pi_2^{-1}(x))} (\pi_2)_! \id_{\pi_1^{-1}(Y_\mathbf{a})} (x) \\
&= q^{-\dim (\pi_2^{-1}(x))} \chi_q(\pi_2^{-1}(x) \cap \pi_1^{-1}(Y_\mathbf{a})).
\end{align*}
Now, we already know from the above discussion that $\pi_2^{-1}(x)
\cap \pi_1^{-1} (Y_\mathbf{a}) = \emptyset$ if $x \not \in
Y_{\mathbf{a}^+}$.  So
assuming $x = (\mathbf{D},W,t) \in Y_{\mathbf{a}^+}$, let $r=\rank t$.  Then
\begin{align*}
F\id_{Y_\mathbf{a}} (\mathbf{D},W,t) &= q^{-\dim (\pi_2^{-1}(\mathbf{D},W,t))}
\chi_q(\pi_2^{-1}(\mathbf{D},W,t) \cap \pi_1^{-1}(Y_\mathbf{a})) \\
&= q^{-\dim \Proj^{\s{\mathbf{a}^+} - r -1}} \chi_q \left( \Proj^{\s{\mathbf{a}^+} -
    r -1} \right) \\
&= q^{-(\s{\mathbf{a}^+} - r - 1)} \sum_{i=0}^{\s{\mathbf{a}^+} - r - 1} q^{2i}
\\
&= \left[ \s{\mathbf{a}^+} - r \right] \\
&= [(\mbox{\# down arrrows in $\mathbf{a}^+$}) - (\mbox{\# lower curves in
  $\mathbf{a}^+$})] \\
&= [\mbox{\# unmatched down arrows in $\mathbf{a}^+$}].
\end{align*}
Thus,
\begin{equation}
F\id_{Y_\mathbf{a}} = [\mbox{\# unmatched down arrows in $\mathbf{a}^+$}]
\id_{Y_{\mathbf{a}^+}}.
\end{equation}
Now,
\begin{align*}
E\id_{Y_\mathbf{a}}(x) &= q^{-\dim (\pi_1^{-1}(x))} (\pi_1)_! \pi_2^*
\id_{Y_\mathbf{a}}(x) \\
&= q^{-\dim (\pi_1^{-1}(x))} (\pi_1)_! \id_{\pi_2^{-1}(Y_\mathbf{a})} (x) \\
&= q^{-\dim (\pi_1^{-1}(x))} \chi_q(\pi_1^{-1}(x) \cap \pi_2^{-1}
(Y_\mathbf{a})).
\end{align*}
We know that $\pi_1^{-1}(x) \cap \pi_2^{-1} (Y_\mathbf{a}) = \emptyset$ if $x
\not \in Y_{\mathbf{a}^-}$.  So assuming $x = (\mathbf{D},W,t) \in Y_{\mathbf{a}^-}$, let
$r = \rank t$.  Then
\begin{align*}
E\id_{Y_\mathbf{a}}(\mathbf{D},W,t) &= q^{-\dim (\pi_1^{-1}(\mathbf{D},W,t))}
\chi_q(\pi_1^{-1}(\mathbf{D},W,t) \cap \pi_2^{-1}(Y_\mathbf{a})) \\
&= q^{-\dim \Proj^{d - r - \s{\mathbf{a}^-} - 1}} \chi_q
(\Proj^{d - r - \s{\mathbf{a}^-} - 1}) \\
&= [d - r - \s{\mathbf{a}^-}] \\
&= [d - (\mbox{\# lower curves in $\mathbf{a}^-$}) - (\mbox{\# down arrows
  in $\mathbf{a}^-$})] \\
&= [(\mbox{\# up arrows in $\mathbf{a}^-$}) - (\mbox{\# lower curves in
  $\mathbf{a}^-$})] \\
&= [\mbox{\# unmatched up arrows in $\mathbf{a}^-$}].
\end{align*}
Thus,
\begin{equation}
E\id_{Y_\mathbf{a}} = [\mbox{\# unmatched up arrows in $\mathbf{a}^-$}]
\id_{Y_{\mathbf{a}^-}}.
\end{equation}
Finally, it is easy to see that
\begin{equation}
\begin{split}
K\id_{Y_\mathbf{a}} &= q^{\pm (d - 2\s{\mathbf{a}})} \id_{Y_\mathbf{a}} \\
&= q^{\pm (\mu - 2(\text{\# unmatched down arrows in $\mathbf{a}$}))}
\id_{Y_\mathbf{a}}
\end{split}
\end{equation}
where $\mu$ is the total number of unmatched arrows in $\mathbf{a}$.
Using the fact that $\mu$ is the total number of middle curves of $b$
(and hence the total number of unmatched vertices in any $\mathbf{a}$
such that ${\bar {\mathbf{a}}} = b$),
the second statement of the theorem now follows easily from a
comparison with \eqref{U_action_V_d}.

Since we know from Section~\ref{sec:crossmatch} that the set \cm{\mu}\ is in
one to one correspondence with the set of intertwiners \twine{\mu},
we have that
\[
\mathcal{T}_s(\mathbf{d}) \cong 
\bigoplus_{\mu} \twine{\mu} \otimes V_\mu \cong V_{\mathbf{d}_1} \otimes
\cdots \otimes V_{\mathbf{d}_k}
\]
which proves the first statement of the theorem.
\end{proof}

Now, like the canonical basis, the basis $\B_s$ we have constucted here
is closely related to the irreducible components of \tv.  To see this,
we first need
a proposition.  
Consider the varieties
$Y_\mathbf{a}$ and $Y_b$ defined over \barFq.  To avoid confusion, denote 
these by $Y_\mathbf{a}'$ and $Y_b'$.  Then

\begin{prop}
$\overline{Y_\mathbf{a}'} = \overline{Z_\mathbf{a}'}$.
\end{prop}

\begin{proof}
Since the $Y_\mathbf{a}'$ are smooth and connected, they are irreducible.
Also, from an argument analogous to the one given in the proof of
Proposition \ref{prop:Y_b_union}, we know that
$\sqcup_\mathbf{a} Y_\mathbf{a}' = \tv'$.  Thus, since the cardinality
of the sets 
$\{\overline{Y_\mathbf{a}'}\}$ and $\{\overline{Z_\mathbf{a}'}\}$ are the same,
$\{\overline{Y_\mathbf{a}'}\}$ must be the set of irreducible components of
$\tv'$.  Now,
$Y_\mathbf{a}' \cap Z_\mathbf{a}' = \bigcup_{\mathbf{r}} A_{\mathbf{a}, \mathbf{r},
  \mathbf{n}^{M(\mathbf{d}, \mathbf{a})}}'$.  But, by
Proposition \ref{prop:A_irred_relation_1},
$\overline{A_{\mathbf{a}, \mathbf{r}^{M(\mathbf{d}, \mathbf{a})},
  \mathbf{n}^{M(\mathbf{d}, \mathbf{a})}}'} = \overline{Z_\mathbf{a}'}$.
Therefore we must have $\overline{Y_\mathbf{a}'} = \overline{Z_\mathbf{a}'}$.
\end{proof}

Since $Y_\mathbf{a}$ is precisely the set of \Fq\ points of
$Y_\mathbf{a}'$, we have the following characterization of the basis $\B_s$.

\begin{theorem}
\label{thm:char-B_s}
The elements $\id_{Y_\mathbf{a}}$ of the basis $\B_s$ are the
unique elements of $\mathcal{T}_s(\mathbf{d})$ equal to one
on the dense points of the irreducible component
$\overline{Z_\mathbf{a}}$ of \tv\ with support contained in this
irreducible component.
\end{theorem}
So, like
the elements of $\B_c$, the elements of
$\B_s$ are equal to a non-zero constant on the set of dense points of an
irreducible component of \tv\ with supports contained in distinct irreducible
components.  However,
unlike $\B_c$, the elements of $\B_s$ have disjoint supports.


\subsection{The Multiplicity Variety $\mathfrak{S}(\mathbf{d})$}

We briefly describe here the relation between $\B_I$ and $\B_s$ and the
multiplicity variety \cite{malkin}.  Let $\mathbf{d} \in (\Z_{\ge
  0})^k$ and let $D$ be a $\s{\mathbf{d}}$-dimensional \barFq\ vector
space.  The \emph{multiplicity variety} is the variety (defined over \barFq)
\[
\mathfrak{S}(\mathbf{d})' = \{(\mathbf{D},t)\, |\, (\mathbf{D},W,t) \in
\tv' \text{ for some } W \subset D\}.
\]
Define the projection $\pi: \mathfrak{T}(\mathbf{d})' \to
\mathfrak{S}(\mathbf{d})'$
by $\pi(\mathbf{D},W,t) = (\mathbf{D},t)$.  It follows
easily from the above results that the irreducible components of
$\mathfrak{S}(\mathbf{d})'$ are given by the closures of the sets
\[
\mathcal{Y}_b' = \{(\mathbf{D},t)\, |\, \alpha(t, \mathbf{D}) =
\mathbf{n}^b\},\ b \in \lcm,
\]
and that these irreducible components are in one to one correspondence
with the irreducible modules in the direct sum decomposition of
$V_{\mathbf{d}_1} \otimes \dots \otimes V_{\mathbf{d}_k}$.  Then $Y_b'
= \pi_{-1}(\mathcal{Y}_b')$ and $\{Y_\mathbf{a}\, |\, {\bar {\mathbf{a}}}
= b\}$ yields
a decomposition of the \Fq\ points of the fiber of $\pi|_{Y_b'}$
isomorphic to the
decomposition of $\qreptr{d}$ into the
subsets $\qreptcr{w}{d}$ where $r$ is the number of lower curves in
$b$.  Thus the bases $\B_I$ and $\B_s$ have natural geometric
interpretations in terms of the multiplicity variety and the
projection $\pi$.


\subsection{The Action of the Intertwiners on $\mathcal{T}_s(\mathbf{d})$}
\label{sec:intertwiner-action}
We will now determine how our intertwiners act on the space
$\mathcal{T}_s(\mathbf{d})$.
For $\mathbf{a} \in (Z_{\ge 0})^k$, let $\mathbf{a}^j = \mathbf{a}^{(1,j)}$.
We will need the following two technical
lemmas.

\begin{lemma}
\label{lemma:flags}
If $\mathbf{D} = (0=\mathbf{D}_0 \subset \mathbf{D}_1 \subset
\mathbf{D}_2 \subset \dots \subset \mathbf{D}_k = D)$ is a
flag with $\mathbf{d} = \alpha(D,\mathbf{D})$ and $\mathbf{a} \in
(\Z_{\ge 0})^k$ with $\mathbf{a}_i \le \mathbf{d}_i$, then
\[
\chi_q(\{W\, |\, W \subset D,\, \alpha(W,\mathbf{D}) = \mathbf{a} \})
= c_{\mathbf{d},\mathbf{a}}
\stackrel{\text{def}}{=} \sum_{\mathbf{b} \in C_\mathbf{a}} q^{2
  \sum_{1 \le j < i \le d} \mathbf{b}_i (1-\mathbf{b}_j)}
\]
where
\[
C_\mathbf{a} = \{ \mathbf{b} \in (\Z_{\ge 0})^d |\ \mathbf{b}_i \in \{0,1\} \,
  \forall \, i,\, \mathbf{b}^{(\mathbf{d}_{j-1}
  + 1, \mathbf{d}_j)} = \mathbf{a}_j \}
\]
and we set $\mathbf{d}_0 = 0$.
\end{lemma}

\begin{proof}
Complete $\mathbf{D}$ to a flag $\mathbf{F}= (0 \subset \mathbf{F}_1
\subset \mathbf{F}_2 \subset \cdots
\subset \mathbf{F}_d = D)$ such that $\dim \mathbf{F}_i = i$ and
$\mathbf{F}_{\mathbf{d}^i}=\mathbf{D}_i$ where $d = \s{\mathbf{d}}$.
This gives a
decomposition of $\gr{\s{\mathbf{a}}}{d}$ into
cells, each isomorphic to $(\Fq)^j$ for some $j$.  The cells are given by
$\{W\, |\, W \subset D,\, \alpha(W,\mathbf{F}) = \mathbf{b} \}$ for a
fixed $\mathbf{b}$. The number of points in such a cell is equal to
\[
q^{2\sum_{1 \le j < i \le d} \mathbf{b}_i (1-\mathbf{b}_j)}.
\]
Our variety is the union of those cells such that
$\mathbf{b}^{(\mathbf{d}_{j-1}+1, \mathbf{d}_j)} = \mathbf{a}_j$.
The result follows.
\end{proof}

Specializing to $q=1$ yields

\begin{lemma}
\label{lemma:flags_q=1}
\begin{align*}
c_{\mathbf{d},\mathbf{a}}
\vert_{q=1} = \prod_{i=1}^k \binom{\mathbf{d}_i}{\mathbf{a}_i}.
\end{align*}
\end{lemma}

\begin{proof}
This follows immediately from Lemma~\ref{lemma:flags} since
$\mathbf{b}_i \in \{0,1\}$ for each cell.
\end{proof}

\begin{theorem}
\label{thm:intertwiner_action_T_s}
The set $\B_I$ acting on $\mathcal{T}_s(\mathbf{d})$
spans the space of
intertwiners \linebreak
$\bigoplus_\mu \twine{\mu}$.  In particular, for $b \in \cm{\mu}$,
$T_{Y_b}$ acts on the basis $\B_s$ of $\mathcal{T}_s(\mathbf{d})$ as
\[
T_{Y_b} \id_{Y_\mathbf{a}} = 
  \begin{cases}
    c_b \id_{\qreptcr{\s{\mathbf{a}}}{d}} \in \mathcal{M}^r(d)
    \cong V_{d - 2r} = V_\mu& \text{if $\bar{\mathbf{a}} = b$}, \\
    0& \text{if $\bar{\mathbf{a}} \ne b$}
  \end{cases}
\]
where $r$ is the
number of lower curves in $b$, $d = \s{\mathbf{d}}$ and $c_b$ is
non-zero constant.
\end{theorem}

\begin{proof}
Recall that $T_Y = p_!
R_Y$.  It is
obvious from the fact that $Y_b = \bigcup_{\mathbf{a} : \bar{\mathbf{a}} = b}
Y_\mathbf{a}$ that
\[
R_{Y_b} \id_{Y_\mathbf{a}} = \id_{Y_b} \id_{Y_\mathbf{a}} = \left\{
  \begin{array}{ll} \id_{Y_\mathbf{a}} & \mbox{if $\bar{\mathbf{a}} = b$} \\ 0 &
    \mbox{if $\bar{\mathbf{a}} \ne b$} \end{array} \right. .
\]
So we need only determine $p_! \id_{Y_\mathbf{a}}$ for $\bar{\mathbf{a}} =
b$.  Now, for $x=(W^x,t^x)
\in \qrept$,
\[
p_! \id_{Y_\mathbf{a}} (x) = \chi_q(p^{-1} (x) \cap Y_\mathbf{a}).
\]
Recall that $p$ is the map $(\mathbf{D},W,t) \mapsto (W,t)$ and
\begin{align*}
Y_\mathbf{a} &= \{(\mathbf{D},W,t)\in \tv \, |\, \alpha(t,\mathbf{D}) =
\mathbf{n}^b,\, \dim W = \s{\mathbf{a}} \}. \\
\end{align*}
Thus,
\begin{equation}
\label{eq:c_b-variety}
p^{-1}(x) \cap Y_\mathbf{a} \cong \{\mathbf{D}\, |\, \dim
(\mathbf{D}_i/ \mathbf{D}_{i-1}) =
\mathbf{d}_i,\, t^x(\mathbf{D}_i) \subset \mathbf{D}_{i-1},\,
\alpha(t^x,\mathbf{D}) = \mathbf{n}^b \}
\end{equation}
if $\dim W^x = \s{\mathbf{a}}$
and $p^{-1}(x) \cap Y_\mathbf{a} = \emptyset$ otherwise.  Note that this
variety depends only on the dimension of the kernel of $t^x$ (or
equivalently, the rank of $t^x$) and the
dimension of $W^x$.  The variety is empty unless $r=\rank t^x$ is
equal to the number of lower curves in $\mathbf{a}$.  Thus, $T_{Y_b}
\id_{Y_\mathbf{a}}$ is a constant function on
$\qreptcr{\s{\mathbf{a}}}{d}$.   Moreover, this constant $c_b$, equal
to the number of points in the variety in \eqref{eq:c_b-variety},
depends only on ${\bar{\mathbf{a}}} = b$ and not on the orientation of
$\mathbf{a}$.
As long as $c_b$ is non-zero, we know that $T_{Y_b}$ is a
non-zero intertwiner.  Moreover, it is obvious that if all the $Y_b$
are non-zero then the intertwiners $T_{Y_b}$ are linearly
independent.

To show that $c_b \ne 0$ it suffices to show that its evaluation at
$q=1$ is non-zero.
The variety \eqref{eq:c_b-variety} consists
of all $t_x$-stable flags $\mathbf{D} = (0 \subset \mathbf{D}_1
\subset \cdots \subset \mathbf{D}_k = D)$
such that $\dim \mathbf{D}_i = \mathbf{d}^i$ and the intersection of
$\mathbf{D}_i$ with $\ker t^x$ is a space of dimension $(\mathbf{n}^b)^j
= \sum_{i=1}^j \mathbf{n}^b_i$.  There is
only one choice for $\mathbf{D}_k$, namely $D$.  Assume we have picked
$\mathbf{D}_{j+1}$.  $\mathbf{D}_j$ can be any subspace of dimension $\mathbf{d}^i$ such that
\begin{align*}
& t^x(\mathbf{D}_{j+1}) \subset \mathbf{D}_j \subset \mathbf{D}_{j+1} \\
& \mbox{and } \dim (\mathbf{D}_j \cap \ker t^x) = (\mathbf{n}^b)^j.
\end{align*}
Note that since $\dim (\mathbf{D}_{j+1} \cap \ker t^x) = (\mathbf{n}^b)^{j+1}$ and
  $\dim \mathbf{D}_{j+1} = \mathbf{d}^{j+1}$, we have that $\dim
  t^x(\mathbf{D}_{j+1}) = \mathbf{d}^{j+1} - (\mathbf{n}^b)^{j+1}$. Also, since
  $(t^x)^2=0$, $t^x(\mathbf{D}_{j+1})\subset \ker t^x$.  Passing to
  the quotient by $t^x(\mathbf{D}_{j+1})$ and denoting this by a bar,
  we see that
  picking a subspace $\mathbf{D}_j$ subject to the above conditions is
  equivalent to picking a subspace $\overline{\mathbf{D}_j}$ of
$\overline{\mathbf{D}_{j+1}}$ of dimension $\mathbf{d}^j - (\mathbf{d}^{j+1} -
  (\mathbf{n}^b)^{j+1})$ such that
\[
\dim (\overline{\mathbf{D}_j} \cap \overline{\ker t^x}) = (\mathbf{n}^b)^j -
(\mathbf{d}^{j+1} - (\mathbf{n}^b)^{j+1}).
\]
Since $\dim \overline{\mathbf{D}_{j+1}} = \mathbf{d}^{j+1} - (\mathbf{d}^{j+1} -
(\mathbf{n}^b)^{j+1}) =  (\mathbf{n}^b)^{j+1}$ and $\dim
\overline{\mathbf{D}_{j+1}}\cap \overline{\ker t^x} = \dim \overline{\mathbf{D}_{j+1}
  \cap \ker t^x} = (\mathbf{n}^b)^{j+1} -
(\mathbf{d}^{j+1} - (\mathbf{n}^b)^{j+1}) = 2(\mathbf{n}^b)^{j+1} -
\mathbf{d}^{j+1}$ we see by Lemma~\ref{lemma:flags_q=1} that the value
of $\chi_q$ of the variety of such spaces evaluated at $q=1$ is
\begin{align*}
& \binom{2(\mathbf{n}^b)^{j+1}  - \mathbf{d}^{j+1}} {(\mathbf{n}^b)^{j+1} +
    (\mathbf{n}^b)^j - \mathbf{d}^{j+1}} \cdot
    \binom{(\mathbf{n}^b)^{j+1} - (2(\mathbf{n}^b)^{j+1} -
    \mathbf{d}^{j+1})} {\mathbf{d}^j - \mathbf{d}^{j+1} +
    (\mathbf{n}^b)^{j+1} - ((\mathbf{n}^b)^{j+1} + (\mathbf{n}^b)^j
    - \mathbf{d}^{j+1})} \\
&= \binom{2(\mathbf{n}^b)^{j+1} - \mathbf{d}^{j+1}} {(\mathbf{n}^b)^{j+1}
    + (\mathbf{n}^b)^j - \mathbf{d}^{j+1}} \cdot
   \binom{\mathbf{d}^{j+1} - (\mathbf{n}^b)^{j+1}} {\mathbf{d}^j -
    (\mathbf{n}^b)^j}.
\end{align*}
This is thus strictly positive provided that
\begin{align}
2(\mathbf{n}^b)^{j+1} - \mathbf{d}^{j+1} &\ge 0 \label{ineq_1}  \\
(\mathbf{n}^b)^{j+1} + (\mathbf{n}^b)^j - \mathbf{d}^{j+1}
&\ge 0 \label{ineq_2} \\
\mathbf{d}^{j+1} - (\mathbf{n}^b)^{j+1} &\ge 0 \label{ineq_3} \\
\mathbf{d}^j - (\mathbf{n}^b)^j &\ge 0 \label{ineq_4} \\
2(\mathbf{n}^b)^{j+1} - \mathbf{d}^{j+1} &\ge (\mathbf{n}^b)^{j+1} +
(\mathbf{n}^b)^j - \mathbf{d}^{j+1} \label{ineq_5} \\
\mathbf{d}^{j+1} - (\mathbf{n}^b)^{j+1} &\ge \mathbf{d}^j -
(\mathbf{n}^b)^j \label{ineq_6}
\end{align}
Now, recall that $\mathbf{n}^b$ is obtained from $\mathbf{a}$ by
forcing all unmatched
arrows to be oriented down.  Also, $\mathbf{d}^j$ is the number of
vertices associated to $V_{\mathbf{d}_1}$ through $V_{\mathbf{d}_j}$ while
$(\mathbf{n}^b)^j$ is number of these vertices with down arrows.
Thus $\mathbf{d}^{j} - (\mathbf{n}^b)^j$ is the number of these
vertices with up arrows.  So \eqref{ineq_3}, \eqref{ineq_4} and \eqref{ineq_6}
are obvious.  \eqref{ineq_5} follows from the simple
fact that $(\mathbf{n}^b)^{j+1} \ge (\mathbf{n}^b)^j$.
\eqref{ineq_1} and \eqref{ineq_2} follow from the fact that each up
arrow is matched to a down arrow to its left since all unmatched
arrows point down and matchings are oriented to the left.

Thus, $\chi_q$ of the variety of choices of $\mathbf{D}_j$ given
$\mathbf{D}_{j+1}$ is independent
of $\mathbf{D}_{j+1}$ (up to isomorphism) and is non-zero. 
Using the fact that the Euler characteristic of a locally trivial fibered space
is equal to the product of the Euler characteristics of the base and
the fiber, we see that the evaluation of $c_b$ at 1 is a product of
positive numbers and is thus
positive.  So $c_b \ne 0$.
\end{proof}


\subsection{The Action of the Intertwiners on $\mathcal{T}_c(\mathbf{d})$}
We now compute the action of our intertwiners on the space
$\mathcal{T}_c(\mathbf{d})$.

Define the coefficients $\kappa^{\mathbf{d},\mathbf{w}}_\mathbf{a}$ by
\[
\diamondsuit \sv{\mathbf{w}} = \sum_\mathbf{a} \kappa^{\mathbf{d},
  \mathbf{w}}_\mathbf{a} \left( \otimes \sv{\mathbf{a}} \right) .
\]
For $b \in \cm{\mu}$, define $\mathbf{l}^b, \mathbf{m}^b \in (\Z_{\ge
  0})^k$ such that $\mathbf{l}^b_i$ is equal to the number of left endpoints
of lower curves of $b$ in the box corresponding to $V_{\mathbf{d}_i}$
and $\mathbf{m}^b_i$ is equal to the number of endpoints of middle
curves of $b$ in the box corresponding to $V_{\mathbf{d}_i}$.

\begin{theorem}
\label{thm:intertwiner_action_B_c}
The set $\B_I$ acting on $\mathcal{T}_c(\mathbf{d})$ spans the space
of intertwiners
\linebreak
$\bigoplus_\mu \twine{\mu}$.  In particular,
if $b \in \cm{\mu}$ is such that $b \le M(\mathbf{d},\mathbf{w})$, then
\[
  T_{Y_b} (g^{\mathbf{d}}_\mathbf{w}) = \omega
  \id_{\qreptcr[\s{\mathbf{l}^b}]{\s{\mathbf{w}}}{d}}
\]
where
\begin{align*}
\omega &= \sum_\mathbf{a} \left( \kappa^{\mathbf{m}^b, \mathbf{w} -
  \mathbf{l}^b}_\mathbf{a} k_{\mathbf{a}+\mathbf{l}^b, \mathbf{l}^b,
  \mathbf{l}^b + \mathbf{m}^b}
  \prod_{i=1}^{k-1} c_{\mathbf{a}_1^i, \mathbf{a}_2^i} \right), \\
\mathbf{a}_1^i &=  \left( (\mathbf{l}^b)^{(i,k)},\, \mathbf{a}^{(i,k)},\,
  (\mathbf{m}^b - \mathbf{a})^{(i,k)},\,
  (\mathbf{d} - \mathbf{m}^b - 2\mathbf{l}^b)^{(i,k)} \right), \\
\mathbf{a}_2^i &= \left( \mathbf{l}^b_i,\, \mathbf{a}_i,\, \mathbf{m}^b_i -
  \mathbf{a}_i,\, \mathbf{d}_i - \mathbf{m}^b_i - \mathbf{l}^b_i \right)
\end{align*}
Otherwise, $T_{Y_b} (g^{\mathbf{d}}_\mathbf{w}) = 0$.
\end{theorem}

\begin{proof}
For a crossingless match $b \in \cm{\mu}$,
$T_{Y_b} (g^{\mathbf{d}}_{\mathbf{w}}) = p_! R_{Y_b}
(g^{\mathbf{d}}_\mathbf{w})$ and
\begin{equation*}
R_{Y_b} (g^{\mathbf{d}}_{\mathbf{w}}) = \sum_{S \le M(\mathbf{d},
  \mathbf{w})} R_{Y_b} \left(
  \eta_{\mathbf{r}^S, \mathbf{n}^S} \right)^{-1} \left( \diamondsuit \sv[\mathbf{n}^S -
  \mathbf{r}^S]{\mathbf{w} - \mathbf{r}^S} \right) .
\end{equation*}
This is equal to zero unless $\bar S = b$ for some $S \
\le M(\mathbf{d}, \mathbf{w})$ (that is, the set of lower curves of $b$
is a subset of the set of lower curves of $M(\mathbf{d}, \mathbf{w})$).
If this is the case, then
\begin{equation*}
R_{Y_b} (g^{\mathbf{d}}_{\mathbf{w}}) = \left(
  \eta_{\mathbf{r}^S, \mathbf{n}^S} \right)^{-1} \left( \diamondsuit
  \sv[\mathbf{n}^S -
  \mathbf{r}^S]{\mathbf{w} - \mathbf{r}^S} \right)
\end{equation*}
for the particular $S \le M(\mathbf{d}, \mathbf{w})$ such that $\bar S
= b$.  Then
$\mathbf{n}^S = \mathbf{l}^b + \mathbf{m}^b$ and
$\mathbf{r}^S = \mathbf{l}^b$.  So
\begin{align*}
R_{Y_b} (g^{\mathbf{d}}_{\mathbf{w}}) &= \left(
  \eta_{\mathbf{l}^b, \mathbf{l}^b+\mathbf{m}^b} \right)^{-1} \left( \diamondsuit \sv[\mathbf{m}^b]{\mathbf{w} -
  \mathbf{l}^b} \right) \\
&= \left( \eta_{\mathbf{l}^b, \mathbf{l}^b+\mathbf{m}^b} \right)^{-1}
\left( \sum_{\mathbf{a}} \kappa^{\mathbf{m}^b, \mathbf{w} -
  \mathbf{l}^b}_{\mathbf{a}} \left( \otimes \sv[\mathbf{m}^b]{\mathbf{a}}
  \right) \right) \\
&= \sum_\mathbf{a} \kappa^{\mathbf{m}^b, \mathbf{w} - \mathbf{l}^b}_\mathbf{a}
f_{\mathbf{a}+\mathbf{l}^b, \mathbf{l}^b, \mathbf{l}^b + \mathbf{m}^b} \\
&= \sum_\mathbf{a} \kappa^{\mathbf{m}^b, \mathbf{w} - \mathbf{l}^b}_\mathbf{a} k_{\mathbf{a}+\mathbf{l}^b, \mathbf{l}^b, \mathbf{l}^b + \mathbf{m}^b}
\id_{A_{\mathbf{a}+\mathbf{l}^b, \mathbf{l}^b, \mathbf{l}^b + \mathbf{m}^b}} .
\end{align*}

Let $(W,t) \in \qrept[d]$.  Then if the set of lower curves of
$b$ is a subset of the set of lower curves of $M(\mathbf{d}, \mathbf{w})$,
\begin{align}
T_{Y_b} (g^{\mathbf{d}}_\mathbf{w}) (W,t) &= \sum_\mathbf{a}
  \kappa^{\mathbf{m}^b, \mathbf{w}
  - \mathbf{l}^b}_\mathbf{a} k_{\mathbf{a}+\mathbf{l}^b, \mathbf{l}^b,
  \mathbf{l}^b + \mathbf{m}^b}
  p_! \id_{A_{\mathbf{a} + \mathbf{l}^b, \mathbf{l}^b, \mathbf{l}^b + \mathbf{m}^b}} (W,t) \nonumber \\
\label{eq:intertwiner_can_coeff}
\begin{split}
&= \sum_\mathbf{a} \kappa^{\mathbf{m}^b, \mathbf{w} - \mathbf{l}^b}_\mathbf{a} k_{\mathbf{a}+\mathbf{l}^b, \mathbf{l}^b, \mathbf{l}^b +
  \mathbf{m}^b} \chi_q (\{\mathbf{D} \, |\,
  \alpha(D,\mathbf{D}) = \mathbf{d},\, t(\mathbf{D}_i) \subset
  \mathbf{D}_{i-1}, \\
&\qquad
  \alpha(\im t,\mathbf{D}) = \mathbf{l}^b,\,
  \alpha(W,\mathbf{D}) = \mathbf{l}^b + \mathbf{a},\,
  \alpha(\ker t, \mathbf{D}) = \mathbf{l}^b + \mathbf{m}^b\}).
\end{split}
\end{align}
We see from Proposition~\ref{prop:can_basis} that $\kappa^{\mathbf{m}^b,
  \mathbf{w} - \mathbf{l}^b}_\mathbf{a} = 0$ unless $\s{\mathbf{a}} =
\s{\mathbf{w} - \mathbf{l}^b} = \s{\mathbf{w}} - \s{\mathbf{l}^b}$.  Therefore, since
$\s{\alpha(W,\mathbf{D})} = \dim W$, \eqref{eq:intertwiner_can_coeff} is zero
unless $\dim W = \s{\mathbf{w}}$.  Similarly, it is zero unless $\rank t =
\dim (\im t) = \s{\mathbf{l}^b}$.  If these conditions are satisfied,
\eqref{eq:intertwiner_can_coeff} is independent of $W$ and $t$.
We can then evaluate $\omega$, the value of the expression in
\eqref{eq:intertwiner_can_coeff}, using Lemma \ref{lemma:flags} and
the fact that the Euler characteristic of a locally trivial
fibered space is the product of the Euler characteristics of
the base and the fiber.  There is only one possible choice for
$\mathbf{D}_0$, namely 0.  Assume we have picked $\mathbf{D}_{i-1}$.
Then $\mathbf{D}_i$ must satisfy the following conditions:
\begin{enumerate}
\item $t(\mathbf{D}_i) \subset \mathbf{D}_{i-1}$ or, equivalently, 
  $\mathbf{D}_i \subset t^{-1}(\mathbf{D}_{i-1})$
\item $\mathbf{D}_i \supset \mathbf{D}_{i-1}$, $\dim \mathbf{D}_i
  = \mathbf{d}^{(1,i)}$
\item $\dim (\mathbf{D}_i \cap \im t) = (\mathbf{l}^b)^{(1,i)}$
\item $\dim (\mathbf{D}_i \cap W) = (\mathbf{l}^b + \mathbf{a})^{(1,i)}$
\item $\dim (\mathbf{D}_i \cap \ker t) = (\mathbf{l}^b +
  \mathbf{m}^b)^{(1,i)}$ .
\end{enumerate}
Pass to the quotient by $\mathbf{D}_{i-1}$ and denote this by a
bar.
Let $\mathbf{F}$ be the flag
\[
\mathbf{F} = (\mathbf{F}_0 = 0 \subset \mathbf{F}_1 = {\overline {\im
    t}} \subset \mathbf{F}_2 = {\overline{W}} \subset \mathbf{F}_3 =
    {\overline{\ker t}}, \mathbf{F}_4 = {\overline {t^{-1}(\mathbf{D}_{i-1})}}).
\]
Then the above conditions are equivalent to picking ${\overline
  {\mathbf{D}_i}} \subset
{\overline {t^{-1}(\mathbf{D}_{i-1})}}$ such that
\[
\alpha(\mathbf{D}_i,\mathbf{F}) = (\mathbf{l}^b_i,
  \mathbf{a}_i, \mathbf{m}^b_i - \mathbf{a}_i, \mathbf{d}_i -
  \mathbf{m}^b_i - \mathbf{l}^b_i ) .
\]
Since
\begin{align*}
\dim {\overline {\im t}} &= (\mathbf{l}^b)^{(i,k)} ,\\
\dim {\overline {W}} &= (\mathbf{l}^b + \mathbf{a})^{(i,k)} ,\\
\dim {\overline {\ker t}} &= (\mathbf{l}^b + \mathbf{m}^b)^{(i,k)} ,\\
\text{and } \qquad \dim {\overline {t^{-1}(\mathbf{D}_{i-1})}} &= \dim
t^{-1}(\mathbf{D}_{i-1})
- \dim \mathbf{D}_{i-1} \\
&= \dim (\im t \cap D_{i-1}) + \dim (\ker t) - \dim \mathbf{D}_{i-1}
\\
&= (\mathbf{l}^b)^{(1,i-1)} + (d- \s{\mathbf{l}^b}) - \mathbf{d}^{(1,i-1)} \\
&= (\mathbf{d} - \mathbf{l}^b)^{(i,k)} ,
\end{align*}
the form of the action of the elements of $\B_I$ follows.

It remains to show that the set $\B_I$ spans the space of intertwiners
$\bigoplus_\mu \twine{\mu}$.
Since it follows from Theorem~\ref{thm:intertwiner_action_T_s} that
the cardinality of $\B_I$ is equal to the dimension of $\bigoplus_\mu
\twine{\mu}$, it suffices to show the linear independence of the set
$\B_I$.  Assume that, acting on the space $\mathcal{T}_s(\mathbf{d})$,
\begin{equation}
\label{eq:lin_ind}
\sum_{i} a_i T_{Y_{b_i}} = 0, \quad a_i \ne 0 \ \forall \ i.
\end{equation}
Since the image of $T_{Y_{b_i}}$ is contained in
$\mathcal{M}^{\s{l_{b_i}}}(d)$ by the above results, we may assume that
$\s{l_{b_i}} = \s{l_{b_j}}$ for all $i$ and $j$.  Fix an $i$ and consider a
$\mathbf{w}$ such that $M(\mathbf{d},\mathbf{w}) = b_i$.  All $T_{Y_{b_j}}$,
$j \ne i$, act by zero on $g^{\mathbf{d}}_{\mathbf{w}}$ by the above (since $\s{l_{b_i}} =
\s{l_{b_j}}$, we cannot have $b_j \le b_i =
M(\mathbf{d},\mathbf{w})$).
Also, $T_{Y_{b_i}} \ne 0$ by the above.
Thus $a_i = 0$ which is a contradiction.  Thus the theorem is proved.
\end{proof}

\subsection{An Isomorphism of $\mathcal{T}_c(\mathbf{d})$ with $\mathcal{T}_s(\mathbf{d})$}

For $(\mathbf{D},W,t) \in \tv$, let
\[
B_{\mathbf{D},W,t} = \{(\mathbf{D}',W',t')\, |\, W'=W,\, t'=t,\,
\alpha(t,\mathbf{D}') = \alpha(t,\mathbf{D})\}.
\]
For $f \in \mathcal{T}(\mathbf{d})$ let
\[
\chi_q(f) = \sum_{x \in \tv} f(x).
\]
Let $\xi : \mathcal{T}_c(\mathbf{d}) \to \mathcal{T}_s(\mathbf{d})$ be
the map given by
\[
\xi(f)(\mathbf{D},W,t) = \chi_q (R_{B_{\mathbf{D},W,t}} f).
\]
The fact that the image of $\xi$ is contained in
$\mathcal{T}_s(\mathbf{d})$ follows from the fact that, up to
isomorphism,
$B_{\mathbf{D},W,t}$ depends only on $\alpha(t, \mathbf{D})$ and $\dim
W$.

\begin{prop}
$\xi$ is an \U-module isomorphism.
\end{prop}

\begin{proof}
This follows easily from Theorems \ref{thm:intertwiner_action_T_s} and
\ref{thm:intertwiner_action_B_c} since
\[
\xi = \sum_b \frac{1}{c_b} (T_{Y_b}|_{Y_b})^{-1}
\circ T_{Y_b}.
\]
\end{proof}


%
\bibliographystyle{amsplain}
\bibliography{tensor}
\end{document}